\documentclass[12pt,francais,draft]{smfart}

\usepackage{smfthm,smfenum}
\usepackage[T1]{fontenc}
\usepackage[english,frenchb]{babel}
\addto\extrasfrenchb{%
  \catcode`\;=13
  \catcode`\:=13
  \catcode`\?=13
  \catcode`\!=13
}

\addto\noextrasfrenchb{%
  \catcode`\;=12
  \catcode`\:=12
  \catcode`\?=12
  \catcode`\!=12
}

\usepackage[latin1]{inputenc}
\usepackage{amssymb,mathrsfs}
\let\mathcal\mathscr
\usepackage[matrix,arrow,curve,cmtip]{xy}
\let\oldxymatrixcompile\xymatrixcompile
\def\xymatrixcompile{\selectlanguage{english}\oldxymatrixcompile}
\CompileMatrices

\allowdisplaybreaks[2]

 \calclayout

\let\ldots\dots
\let\cdots\dots
\makeatletter

 \numberwithin{paragraph}{subsection}
 \def\theparagraph {\thesubsection.\arabic{paragraph}}
 \let\c@equation\c@paragraph
 \let\cl@equation\cl@paragraph
 \def\theequation {\thesubsection.\arabic{equation}}
 \let\subsection\Subsection
 \NumberTheoremsAs{paragraph}

 \def\th@plain{%
  \let\thm@indent\noindent
  \thm@headfont{\normalfont\sc}%
  \thm@notefont{\normalfont}
  \thm@preskip.5\linespacing
  \thm@postskip\thm@preskip
  \ifsmf@skippt
      \thm@headpunct{}\let\thmheadnl\@thmheadnl
      \global\smf@skipptfalse
  \else
       \thm@headpunct{\pointrait}
  \fi
  \itshape }

 \def\th@definition{%
  \th@plain
  \thm@headfont{\normalfont\itshape}\normalfont}

 \def\th@remark{%
   \th@plain
   \thm@headfont{\normalfont\itshape}\normalfont}
{\newcount\@hour\newcount\@minute
  \newcount\acl@temp\@hour=\time\@minute=\time
  \divide\@hour by 60\acl@temp=\@hour\multiply\acl@temp by60\relax
  \advance\@minute by-\acl@temp
  \global\edef\clocktime{\the\@hour:\ifnum\@minute<10 0\fi\the\@minute}}

\def\bibliosection{\@startsection{section}{1}%
  \z@{2\linespacing\@plus\linespacing}{.5\linespacing}%
  {\normalfont\bfseries\centering}}

\renewenvironment{thebibliography}[1]{%
  \@xp\bibliosection\@xp*\@xp{\refname}%
  \def\baselinestretch{1}\normalfont\small\labelsep .5em\relax
  \renewcommand\theenumiv{\arabic{enumiv}}\let\p@enumiv\@empty
  \list{\@biblabel{\theenumiv}}{\settowidth\labelwidth{\@biblabel{#1}}%
    \leftmargin\labelwidth \advance\leftmargin\labelsep
    \usecounter{enumiv}}%
    \itemsep 0.1\baselineskip plus0.1\baselineskip
        minus0.1\baselineskip
    \itemindent 0pt
  \sloppy \clubpenalty\@M \widowpenalty\clubpenalty
  \sfcode`\.=\@m
}{%
  \def\@noitemerr{\@latex@warning{Empty `thebibliography' environment}}%
  \endlist
}
\stepcounter{tocdepth}

     \def\l@part{\@tocline{-1}{0pt}{0pt}{}{\bfseries}}
     \def\l@section{\@tocline{1}{0pt}{0pt}{}{}}
     \def\l@subsection{\@tocline{2}{0pt}{12pt}{}{}}
     \def\l@subsubsection{\@tocline{3}{0pt}{}{}{}}
     \def\l@paragraph{\@tocline{5}{0pt}{}{}{}}

\def\baselinestretch{1.1}
\makeatother

\def\V{\mathbf{V}}
\def\P{\mathbf{P}}
\def\A{\mathbf{A}}
\def\R{\mathbf{R}}
\def\Q{\mathbf{Q}}
\def\N{\mathbf{N}}
\def\C{\mathbf{C}}
\def\Z{\mathbf{Z}}
\let\sqtimes\boxtimes
\let\bar\overline
\let\hat\widehat
\let\tilde\widetilde
\let\ra\rightarrow
\let\phi\varphi
\let\eps\varepsilon

\def\Sym{\operatorname{Sym}\nolimits}
\def\Spec{\operatorname{Spec}\nolimits}
\def\Proj{\operatorname{Proj}\nolimits}

\def\Pic{\operatorname{Pic}\nolimits}
\def\res{\operatorname{res}\nolimits}
\def\rg{\operatorname{rang}\nolimits}

\def\hatPic{\mathop{%
  \smash{\widehat{\operatorname{Pic}}}\vphantom{\operatorname{Pic}}}\nolimits}

\def\vol{\operatorname{vol}}

\let\Fourier\check
\def\eff{\mathrm{eff}}
\def\id{\mathrm{id}}
\def\Res{\operatorname{Res}\nolimits}
\def\Tube {{\mathsf {T}}}

\def\gm{\mathbf{G}_m}
\def\ga{\mathbf{G}_a}

\def\Hom{\operatorname{Hom}}

\def\norm#1{\left\|{#1}\right\|}
\def\abs#1{\left|{#1}\right|}

\def\Re{\operatorname{Re}} 
\def\Im{\operatorname{Im}} 

\input yuri.aux
\begin{document}

\title[Fonctions z\^eta des hauteurs des espaces fibr\'es]
{Fonctions z\^eta des hauteurs \\ des espaces fibr\'es}

\author{Antoine Chambert-Loir}
\address{Institut de math\'ematiques de Jussieu\\ Boite 247 \\
4, place Jussieu \\ F-75252 Paris Cedex 05}
\email{chambert@math.jussieu.fr}

\author{Yuri Tschinkel}
\address{Princeton University \\
    Mathematics Department \\ Fine Hall \\ Washington Road \\
   Princeton NJ 08544-1000 \\ USA}
\email{ytschink@math.princeton.edu}

\date{Soumis sur l'arXiv le 2 mars 2000}


\maketitle

{\def\baselinestretch{1}
\tableofcontents
}


\section*{Introduction}

Cet article est le deuxi\`eme d'une s\'erie consacr\'ee
\`a l'\'etude des hauteurs sur
certaines vari\'et\'es alg\'ebriques sur un corps de nombres,
notamment en ce qui concerne la distribution des points rationnels
de hauteur born\'ee. 

Pr\'ecis\'ement, soient $X$ une vari\'et\'e alg\'ebrique projective lisse
sur un corps de nombres $F$, $\mathcal L$ un fibr\'e en droites
sur $X$ et $H_{\mathcal L}:X(\bar F)\ra\mathbf R_+^*$ une fonction
hauteur (exponentielle) pour $\mathcal L$. Si $U$ est un ouvert
de Zariski de $X$, on cherche \`a estimer le nombre
$$ N_U(\mathcal L,H) = \# \{ x\in U(F) \,;\, H_{\mathcal L}(x)\leq H \}$$
lorsque $H$ tend vers $+\infty$.
L'\'etude de nombreux exemples
a montr\'e que l'on peut s'attendre \`a un \'equivalent
de la forme
\begin{equation} \tag{$*$}\label{eq:maninpeyre}%
 N_U(\mathcal L,H) =
    \Theta(\mathcal L) H^{a(\mathcal L)} (\log H)^{b(\mathcal L)-1} 
    (1+o(1)), \quad H\ra+\infty 
\end{equation}
pour un ouvert $U$ convenable et lorsque par exemple $\mathcal L$ et
$\omega_X^{-1}$ (fibr\'e anticanonique) sont amples.
On a en effet un r\'esultat de ce genre lorsque 
$X$ est une vari\'et\'e de drapeaux~\cite{franke-m-t89},
une intersection compl\`ete lisse de bas degr\'e (m\'ethode du cercle),
une vari\'et\'e torique~\cite{batyrev-t98b},
une vari\'et\'e horosph\'erique~\cite{strauch-t99}, 
une compactification équivariante d'un groupe
vectoriel~\cite{chambert-loir-t2000b}, etc.
On dispose de plus d'une description conjecturale
assez pr\'ecise des constantes
$a(\mathcal L)$ et $b(\mathcal L)$ en termes du c\^one des diviseurs
effectifs~\cite{batyrev-m90}
ainsi que de la constante $\Theta(\mathcal L)$
(\cite{peyre95}, \cite{batyrev-t98}).

En fait, on \'etudie plut\^ot la \emph{fonction z\^eta des hauteurs},
d\'efinie par la s\'erie de Dirichlet
$$ Z_U(\mathcal L,s) =  \sum_{x\in U(F)} H_{\mathcal L}(x)^{-s} $$
\`a laquelle on applique des th\'eor\`emes taub\'eriens standard.
Sur cette s\'erie, on peut se poser les questions suivantes :
domaine de convergence, prolongement m\'eromorphe,
ordre du premier p\^ole, terme principal, sans oublier
la croissance dans les bandes verticales \`a gauche du premier p\^ole.
Cela permet de proposer des conjectures de pr\'ecision variable.

Dans cet article, nous consid\'erons certaines fibrations
localement triviales construites de la fa\c con suivante.
Soient $G$ un groupe alg\'ebrique lin\'eaire sur $F$ agissant sur
une vari\'et\'e projective lisse $X$, $B$ une vari\'et\'e projective lisse
sur $F$ et $T$ un $G$-torseur sur $B$ localement trivial pour
la topologie de Zariski. Ces donn\'ees d\'efinissent une vari\'et\'e
alg\'ebrique projective $Y$ munie d'un morphisme $Y\ra B$ dont les
fibres sont isomorphes \`a $X$.
Le c\oe ur du probl\`eme est de comprendre le comportement de
la fonction hauteur lorsqu'on passe d'une fibre \`a l'autre,
comportement vraiment non trivial bien qu'elles soient
toutes isomorphes.

Dans notre premier article (\emph{Torseurs arithmétiques et espaces
fibrés}, \cite{chambert-loir-t99}), nous avons expos\'e en d\'etail
la construction de hauteurs sur de telles vari\'et\'es.
Dans celui-ci, nous appliquons 
ces consid\'erations g\'en\'erales au cas d'une
fibration en vari\'et\'es toriques provenant d'un torseur sous un
\emph{tore d\'eploy\'e,}
pour l'ouvert $U$ d\'efini par le tore.
Nous avons construit les hauteurs
\`a l'aide d'un prolongement du torseur g\'eom\'etrique
en un torseur arithm\'etique,
ce qui correspond en l'occurence
au choix de m\'etriques hermitiennes sur certains fibr\'es en droites.
\'Ecrivons la fonction z\^eta
comme la somme des fonctions z\^eta des fibres 
$$ Z_U (\mathcal L,s) = \sum_{b\in B(F)}  \sum_{x\in U_b(F)} H_{\mathcal
L}(x)^{-s}
 = \sum_{b\in B(F)} Z_{U_b}(\mathcal L|_{U_b},s). $$
Chaque $U_b$ est isomorphe au tore et on peut exprimer
la fonction z\^eta
des hauteurs de $U_b$ \`a l'aide de la formule de Poisson ad\'elique.
De cette fa\c con, la fonction z\^eta de $U$ appara\^{\i}t
comme une int\'egrale
sur certains caract\`eres du tore ad\'elique de la fonction $L$ d'Arakelov
du torseur arithm\'etique sur $B$.

Ainsi, nous pouvons d\'emontrer des \emph{th\'eor\`emes de mont\'ee} :
supposons que $B$ v\'erifie une conjecture,
alors $Y$  la v\'erifie.
Bien s\^ur, la m\'ethode reprend les outils utilis\'es
dans la d\'emonstration de ces conjectures pour les vari\'et\'es toriques
(\cite{batyrev-t98b,batyrev-t95b,batyrev-t96}).

Par exemple, 
nous d\'emontrons au \S\,\ref{subsec:holomorphie},
sous des hypoth\`eses minimales sur $B$, 
l'holomorphie de la fonction
$Z_U(\mathcal L,s)$ pour $\Re(s)>a(\mathcal L)$ ;
cela implique que pour tout $\eps>0$,
le nombre de points rationnels de hauteur
$H_{\mathcal L}$ inf\'erieure \`a $H$ est $O(H^{a(\mathcal L)+\eps})$.
Ensuite, sous des hypoth\`eses raisonnables concernant $B$,
nous \'etablissons un prolongement m\'eromorphe de cette fonction z\^eta
\`a gauche de $a(\mathcal L)$ et nous d\'emontrons que l'ordre du p\^ole
est inf\'erieur ou \'egal \`a $b(\mathcal L)$ ;
cela pr\'ecise la majoration
du nombre de points en $O(H^{a(\mathcal L)}(\log H)^{b(\mathcal L)-1})$.
Enfin, lorsque $\mathcal L=\omega_{Y}^{-1}$, nous d\'emontrons
que le p\^ole est effectivement d'ordre $b(\mathcal L)$ d'o\`u
une estimation de la forme~\eqref{eq:maninpeyre} et nous
identifions la constante $\Theta(\mathcal L)$,
\'etablissant ainsi la conjecture de Manin raffin\'ee par Peyre.
Pour un fibr\'e en droites quelconque, la preuve de la conjecture de
Batyrev--Manin~\cite{batyrev-m90}
avec son raffinement par Batyrev--Tschinkel~\cite{batyrev-t98}
est ramen\'ee \`a la d\'etermination exacte de l'ordre du p\^ole,
c'est-\`a-dire \`a la non-annulation d'une certaine constante.
Dans le cas des vari\'et\'es toriques
ou des vari\'et\'es horosph\'eriques,
l'utilisation de {\og fibrations $\mathcal L$-primitives\fg}
dans~\cite{batyrev-t96} et~\cite{strauch-t99} a permis d'\'etablir
cette conjecture. Moyennant des hypoth\`eses sur $B$, cette m\'ethode
devrait s'\'etendre au sujet de notre \'etude.

Notre m\'ethode impose de disposer de majorations de la fonction
z\^eta des hauteurs (pour $B$) dans les bandes verticales \`a gauche
du premier p\^ole ; nous avons ainsi
t\^ach\'e d'obtenir de telles majorations pour la vari\'et\'e $Y$. 
Il est en outre bien connu que cela entra\^{\i}ne
un d\'eveloppement asymptotique
assez pr\'ecis pour le nombre de points de hauteur born\'ee,
cf. le th\'eor\`eme taub\'erien donn\'e en appendice.
Quelques cas de vari\'et\'es toriques sur $\Q$ avaient en effet 
attir\'e l'attention des sp\'ecialistes
de th\'eorie analytique des nombres
(voir notamment les articles de \'E.~Fouvry
et R.~de la Bret\`eche dans~\cite{peyre-aster},
ainsi que~\cite{breteche98b}).
Notre m\'ethode \'etablit un tel d\'eveloppement pour les vari\'et\'es
toriques lisses, les vari\'et\'es horosph\'eriques, etc.\ sur tout
corps de nombres.

La d\'emonstration de l'existence d'un prolongement m\'eromorphe
de la fonction z\^eta des hauteurs pour les vari\'et\'es toriques
ou pour les vari\'et\'es horosph\'eriques
faisait intervenir un th\'eor\`eme
technique d'analyse complexe \`a plusieurs variables dont
la d\'emonstration se trouve dans~\cite{batyrev-t98b}, \cite{batyrev-t96}
et~\cite{strauch-t99}.
En vue d'obtenir les majorations exig\'ees dans les bandes verticales,
nous sommes oblig\'es d'en pr\'eciser la preuve ; ceci est
l'objet du \S\,\ref{sec:analyse}. 

Dans les \S\,\ref{sec:torique} et \S\,\ref{sec:fibrations}
se situe l'\'etude de la fonction z\^eta des hauteurs
d'une vari\'et\'e torique et d'une fibration en vari\'et\'es
toriques. Pour les vari\'et\'es toriques, 
nous am\'eliorons le terme d'erreur \`a la suite
de \cite{batyrev-t98b,salberger98,breteche98}.
Le th\'eor\`eme de mont\'ee pour les fibrations
g\'en\'eralise le r\'esultat principal de~\cite{strauch-t99}.


\section*{Notations et conventions}

Si $\mathcal X$ est un sch\'ema,
on note $\Pic(\mathcal X)$ le groupe des classes d'isomorphisme
de faisceaux inversibles sur $\mathcal X$.
Si $\mathcal F$ est un faisceau quasi-coh\'erent sur $\mathcal X$,
on note $\V(\mathcal F)=\Spec\Sym \mathcal F$
et $\P(\mathcal F)=\Proj\Sym \mathcal F$ les fibr\'es vectoriels
et projectifs associ\'es \`a $\mathcal F$.

On note $\hatPic(\mathcal X)$
le groupe des classes d'isomorphisme de fibr\'es en droites hermitiens
sur $\mathcal X$ (c'est-\`a-dire des fibrés en droites munis
 d'une m\'etrique hermitienne
continue sur $\mathcal X(\C)$ et invariante par la conjugaison
complexe)..

Si $\mathcal X$ est un $S$-sch\'ema, et si $\sigma\in S(\C)$,
on d\'esigne par $\mathcal X_\sigma$
le $\C$-sch\'ema $\mathcal X\times_\sigma \C$.
Cette notation servira lorsque $S$ est le spectre d'un localis\'e
de l'anneau des entiers d'un corps de nombres $F$, de sorte
que $\sigma$ n'est autre qu'un plongement de $F$ dans $\C$.

Si $G$ est un sch\'ema en groupes sur $S$,
$X^*(G)$ d\'esigne le groupe des $S$-homomorphismes
$G\ra\gm$ (caract\`eres alg\'ebriques).

Si $\mathcal X/S$ est lisse, le faisceau canonique de $\mathcal X/S$,
not\'e $\omega_{\mathcal X/S}$,
est la puissance ext\'erieure maximale de $\Omega^1_{\mathcal X/S}$.

Enfin,
cet article commence au paragraphe~\ref{sec:analyse}.
Les références aux paragrahes~1 et~2 renvoient
ainsi à l'article précédent~\cite{chambert-loir-t99}.


\section{Fonctions holomorphes dans un tube}
\label{sec:analyse}%

Le but de ce paragraphe est de prouver un th\'eor\`eme
d'analyse sur le prolongement m\'eromorphe de
certaines int\'egrales et leur estimation dans des bandes
verticales. Ce th\'eor\`eme g\'en\'eralise un \'enonc\'e
analogue de~\cite{batyrev-t98b,strauch-t99}. La pr\'esentation
en est un peu diff\'erente et le formalisme que nous introduisons
permet de contr\^oler la croissance des fonctions obtenues.
Ce contr\^ole est n\'ecessaire pour utiliser des th\'eor\`emes
taub\'eriens pr\'ecis et am\'eliorer ainsi le d\'eveloppement
asymptotique du nombre de points rationnels de hauteur born\'ee.

Les r\'esultats de ce paragraphe n'interviennent que
dans la preuve des th\'eor\`emes~\ref{theo:torique} et
\ref{theo:montee}.

\subsection{\'Enonc\'e du th\'eor\`eme}

Soit $V$ un $\R$-espace vectoriel r\'eel de dimension finie
muni d'une mesure de Lebesgue $dv$ et d'une norme $\norm\cdot$.
On dispose alors d'une mesure canonique $dv^*$ sur le
dual $V^*$.
Notons $V_\C=V\otimes_\R \C$ le complexifi\'e de $V$.
On appelle \emph{tube} toute partie connexe de $V_\C$
de la forme $\Omega+iV$ o\`u $\Omega$ est une partie connexe de $V$ ; on
le notera $\Tube (\Omega)$.

Soit enfin $M$ un sous-espace vectoriel de $V$
muni d'une mesure de Lebesgue $dm$.

\begin{defi}\label{defi:classe}%
Une \emph{classe de contr\^ole} $\mathcal D$ est la donn\'ee pour tout
couple $M\subset V$ de $\R$-espaces vectoriels de dimension finie
d'un ensemble $\mathcal D(M,V)$ de fonctions mesurables $\kappa:V\ra\R_+$
dites \emph{$\mathcal D({M,V})$-contr\^olantes}
v\'erifiant les propri\'et\'es
suivantes :
\begin {enumerate} \def\theenumi{\alph{enumi}}\def\labelenumi{(\theenumi)}
\item\label{axiome:hypomonoide}%
  si $\kappa_1$ et $\kappa_2$ sont deux fonctions de $\mathcal D(M,V)$,
  $\lambda_1$ et $\lambda_2$ deux réels positifs,
  et si $\kappa$ est une fonction mesurable $V\ra\R_+$
  telle que $\kappa\leq \lambda_1 \kappa_1 +  \lambda_2\kappa_2$,
  alors $\kappa\in \mathcal D(M,V)$ ;
\item\label{axiome:max}%
  Si $\kappa\in\mathcal D(M,V)$ et si $K$ est un compact de $V$, 
  la fonction $v\mapsto \sup_{u\in K} \kappa(v+u)$
  appartient à $\mathcal D(M,V)$ ;
\item\label{axiome:limite}%
  si $\kappa\in\mathcal D(M,V)$,
  pour tout $v\in M\setminus 0$, $\kappa(tv)$ tend vers~$0$
  lorsque $t$ tend vers $+\infty$ ;
\item\label{axiome:integrale}%
  si $\kappa\in\mathcal D(M,V)$,
  pour tout sous-espace $M_1\subset M$, 
  la fonction  $M_1$-invariante
  \[\kappa_{M_1}:v\mapsto \int_{M_1} \kappa(v+m_1)\, dm_1 \]
  est finie et appartient \`a $\mathcal D(M/M_1,V/M_1)$ ;
\item\label{axiome:projecteur}%
  si $\kappa\in\mathcal D(M,V)$,
  pour tout sous-espace $M_1\subset M$ et tout projecteur $p:V\ra V$
  de noyau $M_1$, 
  la fonction $M_1$-invariante $\kappa\circ p$
  appartient \`a $\mathcal D(M/M_1,V/M_1)$.
\end{enumerate}
\end{defi}

\paragraph{}\label{defi:classemax}%
Il existe une classe de contr\^ole $\mathcal D^{\text{max}}$ contenant
toutes les classes de contr\^oles :
l'ensemble $\mathcal D^{\text{max}}(M,V)$
est d\'efini par r\'ecurrence sur la dimension de $M$ par les trois
conditions~(\ref{axiome:hypomonoide}, \ref{axiome:limite},
\ref{axiome:projecteur})
dans la d\'efinition~\ref{defi:classe}.
La derni\`ere condition est alors automatique.

Dans la suite, on fixe une classe de contr\^ole $\mathcal D$,
et on abr\`ege l'expression \emph{$\mathcal D(M,V)$-contr\^olante}
en \emph{$M$-contr\^olante.}

\begin{defi}\label{defi:controlee}%
Une fonction $f:\Tube(\Omega)\ra \C$
sur un tube est dite \emph{$M$-contr\^ol\'ee}
s'il existe une fonction $M$-contr\^olante $\kappa$ telle que pour tout
compact $K\subset\Tube(\Omega)$, il existe un r\'eel $c(K)$
de sorte que l'in\'egalit\'e
\[ \abs{f(z+iv)}\leq c(K) \kappa(v) \]
soit v\'erifi\'ee pour tout $z\in K$ et tout $v\in V$.
\end{defi}

\paragraph{}
Consid\'erons une fonction sur un tube, $f:\Tube (\Omega)\ra\C$.
Soit $M$ un sous-espace vectoriel de $V$, muni d'une mesure
de Lebesgue $dm$.
On consid\`ere la projection $\pi:V\ra V'=V/M$ et on munit $V'$
de la mesure de Lebesgue quotient.
On pose, quand cela a un sens,
\begin{equation}
\mathcal S_M (f) (z)
    = \frac{1}{(2\pi)^{\dim M}} \int_M f(z+im)\, dm, \qquad z\in\Tube (\Omega).
\end{equation}

\begin{lemm}\label{lemm:integrale}%
Soit $\Omega\subset V$ et $f:\Tube(\Omega)\ra\C$ une fonction holomorphe
$M$-contr\^ol\'ee. Soit $M'$ un sous-espace vectoriel
de $M$ et $\Omega'$ l'image de $\Omega$ par la projection
$V\ra V/M'$.
Alors, l'intégrale qui définit $\mathscr S_{M'}(f)$
converge en tout $z\in\Tube(\Omega)$ et
d\'efinit une fonction holomorphe 
$M/M'$-contr\^ol\'ee sur $\Tube(\Omega')$.
\end{lemm}
\begin{proof}
Comme $f$ est $M$-contrôlée, il existe une fonction $\kappa\in\mathscr
D(M,V)$ et, pour tout compact~$K\subset \Tube(\Omega)$,
un réel $c(K)>0$ de sorte que pour tout $v\in V$ et tout $z\in K$,
on ait $\abs{f(z+iv)}\leq c(K)\kappa(v)$.
La condition~(\ref{defi:classe}, \ref{axiome:integrale})
des classes de contrôles jointe
au théorème de convergence dominée de Lebesgue implique que l'intégrale
qui définit $\mathscr S_{M'}(f)$ converge et
que la somme est une fonction holomorphe sur $\Tube(\Omega)$.
Par construction, cette fonction est $iM'$-invariante. Comme elle est
analytique, elle est donc invariante par $M'$ et définit ainsi une fonction 
holomorphe sur $\Tube(\Omega')$.
De plus, si $\pi$ désigne la projection $V\ra V/M'$,
pour tout $z\in K$  et tout $v\in V$, on a
\[ \abs{\mathscr S_{M'}(f)(\pi(z)+i\pi(v)) }
 \leq c(K) \int_{M'} \kappa(v+m')\, dm'
     = c(K) \kappa'(\pi(v)) \]
où $\kappa'$ appartient par définition à $\mathscr D(M/M',V/M')$.
Tout compact de $\Tube(\Omega')$ étant de la forme $\pi(K)$ pour un compact
$K$ de $\Tube(\Omega)$, le lemme est ainsi démontré.
\end{proof}

\paragraph{Fonction caract\'eristique d'un c\^one}
Soit $\Lambda$ un c\^one convexe poly\'edral ouvert de $V$.
La fonction caract\'eristique de $\Lambda$
est la fonction sur $\Tube(\Lambda)$ d\'efinie par l'intégrale
convergente
\begin{equation}
\mathsf X_{\Lambda} (z)
= \int_{\Lambda^*} e^{-\langle z,v^*\rangle}\, dv^* ,
\end{equation}
o\`u $\Lambda^*\subset V^*$ est le c\^one dual de $\Lambda$,
$V^*$ \'etant muni de la mesure de Lebesgue $dv^*$ duale
de la mesure $dv$.

Si $\Lambda$ est simplicial, c'est-\`a-dire qu'il existe 
$n=\dim V$ formes lin\'eaires ind\'ependantes $\ell_1,\ldots,\ell_n$
telles que $v\in\Lambda$ si et seulement si $\ell_j(v)>0$
pour tout $j$, alors 
\[ \mathsf X_{\Lambda}(z) =
\norm{d\ell_1\wedge \cdots \wedge d\ell_n}
\prod_{j=1}^n \frac{1}{\ell_j(z)}. \]
(On a noté $\norm{d\ell_1\wedge \cdots \wedge d\ell_n}$
le volume du parallélépipède fondamental dans $V^*$ de base 
les $\ell_j$.)
Dans le cas général, toute triangulation de $\Lambda^*$
par des cônes simpliciaux permet d'exprimer $\mathsf X_\Lambda$
sous la forme d'une somme de fractions rationnelles de ce type.
Elle se prolonge ainsi en une fonction rationnelle sur $\Tube(V)$
dont les pôles sont exactement les hyperplans de $V_\C$
définis par les équations des faces de $\Lambda$.
Elle est de plus strictement positive sur $\Lambda$.

Une autre façon de construire un cône est de s'en donner des
générateurs, autrement dit de l'écrire comme \emph{quotient}
d'un cône simplicial. \`A ce titre, on a la proposition suivante.

\begin{prop}\label{prop:X-controle}%
Soit $\Lambda$ un c\^one polyédral convexe ouvert de $V$ 
dont l'adhérence $\bar\Lambda$ ne contient pas de droite.
Soit $M$ un sous-espace vectoriel de $V$ tel que $\bar\Lambda\cap M=\{0\}$. 
On note $\pi$ la projection $V\ra V'=V/M$.

La restriction à $\Tube(\Lambda)$ de la fonction $\mathsf X_\Lambda$ 
est $M$-contrôlée (pour la classe $\mathscr D^{\max}(M,V)$).
L'intégrale qui définit $\mathscr S_M(\mathsf X_\Lambda)$
converge donc absolument et pour tout $z\in\Tube(\Lambda)$, on a
\[ \mathcal S_M(\mathsf X_\Lambda) (z)
                 = \mathsf X_{\Lambda'} (\pi (z)) . \]
\end{prop}
\begin{rema}
Les hypothèses impliquent que $\bar\Lambda'$
ne contient pas de droite. En effet, s'il existait un vecteur non nul
de $\bar\Lambda'\cap -\bar\Lambda'$, il existerait deux vecteurs
$v_1$ et $v_2$ de $\bar\Lambda$ tels que $v_1+v_2\in M$
mais $v_1\not\in M$. Comme $\bar\Lambda\cap M=\{0\}$, $v_1=-v_2$
ce qui contredit l'hypothèse que $\bar\Lambda$ ne contient pas
de droite.
\end{rema}

\begin{proof}
La preuve est une adaptation des paragraphes 7.1 et 7.2
de~\cite{strauch-t99}.
Soit $(e_i)$ une famille minimale de générateurs de $\Lambda$.
Chaque face de $\Lambda^*$ dont la dimension est $\dim V-1$
engendre un sous-espace vectoriel qui est l'orthogonal d'un des $e_i$.

Comme $M\cap\bar\Lambda=\{0\}$,
il existe une forme linéaire $\ell\in V^*$ qui est nulle sur $M$
mais qui n'appartient à aucune face de $\Lambda^*$ ;
posons $H=\ker\ell$.
Soit $H'$ un supplémentaire de $\R\ell$ dans $V^*$.
Si $\phi\in V^*$ et $t\in\R$ sont tels que $\phi+t\ell\in\Lambda^*$,
on doit avoir pour tout générateur $e_j$ de $\Lambda$ l'inégalité
$\phi(e_j)+t\ell(e_j)> 0$, soit (rappelons que $\ell(e_j)$ n'est pas nul),
$t > -\phi(e_j)/\ell(e_j)$ quand $\ell(e_j)>0$
et $t< -\phi(e_j)/\ell(e_j)$ quand $\ell(e_j)<0$.
Soit alors $I(\phi)=\left]h_1(\phi),h_2(\phi)\right[$
l'intervalle de $\R$ défini par ces inégalités.
(Si tous les $\ell(e_j)$ sont positifs, c'est-à-dire $\ell\in\Lambda^*$,
on a $h_1\equiv -\infty$, tandis que s'ils sont tous négatifs,
$h_2\equiv +\infty$.)
Les fonctions $h_1$ et $h_2$ sont linéaires par morceaux
par rapport à un éventail de $H'$ qu'on peut supposer complet et régulier
(voir par exemple~\cite{fulton93} pour la définition,
ou~\cite{batyrev-t95b}).

Alors, si $v\in\Tube(\Lambda)$ et $m\in H$, on a
\begin{align*}
\mathsf X_{\Lambda}(v+im)
&= \int_{V^*} \mathbf 1_{\Lambda^*}(\phi)
     e^{-\langle v+im,\phi\rangle}\,d\phi \\
&= \int_{H'} \int_\R \mathbf 1_{\Lambda^*}(\phi+t\ell)
     e^{-\langle v+im,\phi\rangle} e^{-t\langle v,\ell\rangle}\, dt\, d\phi \\
&= \int_{H'} \int_{h_1(\phi)}^{h_2(\phi)}
     e^{-\langle v+im,\phi\rangle} e^{-t\langle v,\ell\rangle}\, dt\, d\phi \\
&= \int_{H'} e^{-\langle v+im,\phi\rangle}
    \frac{e^{- h_1(\phi)\langle v,\ell\rangle}-e^{-h_2(\phi)\langle
v,\ell\rangle}}{\langle v,\ell\rangle} \, d\phi
\end{align*}
de sorte que la fonction $H\ra\C$ telle que $m\mapsto \mathsf X_\Lambda(v+im)$
est (à une constante multiplicative près)
la différence des transformées de Fourier des fonctions
\[ H'\ra\C, \quad \phi\mapsto e^{-\langle v,\phi+h_j(\phi)\ell\rangle } \]
pour $j=1$ et $2$.

Comme $v\in\Tube(\Lambda)$ et $\phi+h_j(\phi)\ell$ appartient
au bord de $\Lambda^*$, $\langle v,\phi+h_j(\phi)\ell\rangle$
est de partie réelle strictement positive, à moins que $\phi=0$.
Soit $K$ un compact de $\Tube(\Lambda)$.
Il résulte alors des estimations des transformées de Fourier
de fonctions linéaires par morceaux et positives
(voir~\cite{batyrev-t95b}, proposition 2.3.2, p.~614, et aussi
\emph{infra}, prop.~\ref{prop:fourierinfini})
une majoration de  la fonction
\[  f_{\Lambda,K} (m) := \sum_{v\in K} \abs{\mathsf X_\Lambda(v+im)} \]
de la forme
\[ f_{\Lambda,K} (m) \leq c(K) \sum_\alpha \prod_{j=1}^{\dim H}
       \frac{1}{(1+\abs{\langle m,\ell_{\alpha,j}\rangle})^{1+1/\dim H}}, \]
où pour tout $\alpha$, la famille $(\ell_{\alpha,j})_j$ est
une base de $H^*$.
D'après le lemme~\ref{lemm:controle-facile} ci-dessous, la fonction
$f_{\Lambda,K}$ appartient à $\mathscr D^{\max}(M,V)$.

La fonction $m\mapsto \mathsf X_\Lambda(v+im)$
est donc absolument intégrable sur $M$.
C'est la transformée de Fourier de la fonction
$\phi\mapsto \mathbf 1_{\Lambda^*}(\phi) e^{-\langle v,\phi\rangle}$ 
dont il est facile de voir qu'elle est intégrable sur tout sous-espace
et donc aussi $M^\perp$.
La formule de Poisson s'applique (après un léger argument de
régularisation) et s'écrit
\[ \int_M \mathsf X_\Lambda(v+im)\, dm
        = (2\pi)^{\dim M}
          \int_{\Lambda^*\cap M^\perp} e^{-\langle v,\phi\rangle} \,d\phi. \]
Or, l'application $V\ra V'$ identifie $(V')^*$ à $M^\perp$,
et $\Lambda^*\cap M^\perp$ à $(\Lambda')^*$.
Ainsi, on obtient
\[ \mathscr S_M (\mathsf X_\Lambda) (v)
       =\int_{(\Lambda')^*} e^{-\langle\pi(v),\phi\rangle}\, d\phi
      = \mathsf X_{\Lambda'} (\pi(v)). \]
\end{proof}

\begin{lemm}\label{lemm:controle-facile}%
Soit $V$ un $\R$-espace vectoriel de dimension $d$,
$(\ell_1,\dots,\ell_d)$ une base de $V^*$ et $f$ la fonction
$v\mapsto \prod_{j=1}^d (1+\abs{\ell_j(v)})^{-1-1/d}$.
Alors, $f\in\mathscr D^{\max}(V,V)$.
\end{lemm}
\begin{proof}
Soit $M$  un sous-espace vectoriel de $V$ de dimension~$m$.
Quitte à réordonner les indices, on peut supposer que
$M$ est l'image d'une application linéaire $\R^m\ra\R^d=V$
de la forme $t=(t_1,\ldots,t_m)\mapsto (t_1,\ldots,t_m,\phi_{m+1}(t),
\ldots,\phi_{d}(t))$.
Si on réalise $V/M$ par son supplémentaire $\{0\}^m\times\R^{d-m}$,
la fonction $f_M:v \mapsto \int_M f(v+m)\, dm$ 
est donnée par l'intégrale
\[ \int_{\R^m} \frac{1}{(1+\abs{t_1})^{1+1/d}}\cdots
\frac{1}{(1+\abs{t_m})^{1+1/d}}
        \prod_{j=m+1}^d \frac{1}{(1+\abs{v_j + \phi_j(t)})^{1+1/d}}\,
dt_1 \cdots dt_m. \]
Elle est dominée par l'intégrale convergente
\[ \int_{\R^m} \frac{1}{(1+\abs{t_1})^{1+1/d}}\cdots 
\frac{1}{(1+\abs{t_m})^{1+1/d}}\, dt_1 \cdots dt_m \]
et le théorème de convergence dominée implique alors que 
pour tout vecteur $v=(0,\ldots,0,v_{m+1},\ldots,v_d)$
distinct de $0$,
\[ \lim_{s\ra+\infty} f_M(sv) = 0. \]
Le lemme est ainsi démontré.
\end{proof}

\begin{defi}\label{defi:H(C)}%
Soient $C$ un ouvert convexe de $V$ ayant $0$ pour point adh\'erent
et $\Lambda$ un c\^one polyédral ouvert contenant $C$.

Soit $\Phi\subset V^*$ une famille de formes linéaires 
deux à deux non proportionnelles définissant les
faces de $\Lambda$.

On note $\mathcal H_M(\Lambda;C)$
l'ensemble des fonctions holomorphes $f:\Tube (C) \ra\C$
telles qu'il existe
un voisinage convexe $B$ de $0$ dans  $V$
de sorte que la fonction $g$ d\'efinie par
\[ g(z) = f(z) \prod_{\phi\in \Phi} \frac{\phi(z)}{1+\phi(z)} \]
admet un prolongement holomorphe $M$-contrôlé dans $\Tube(B)$.
\end{defi}

Par le th\'eor\`eme d'extension de Bochner (voir par
exemple~\cite{narasimhan71}), une telle fonction 
s'\'etend en une fonction holomorphe sur le tube de base l'enveloppe
convexe $C'$ de $B\cup C$. En particulier, il n'aurait pas \'et\'e
restrictif de prendre pour $C$ l'intersection du c\^one $\Lambda$
avec un voisinage convexe de $0$ dans $V$.

On constate aussi que $f$ est nécessairement $M$-contrôlée dans
$\Tube(C)$.
Enfin, il est facile de vérifier que $\mathcal H_M(\Lambda;C)$
ne dépend pas du choix des formes linéaires
qui définissent les faces de $\Lambda$.

\paragraph{}
Si $\Lambda$ est un cône polyédral et si $M$ est un sous-espace vectoriel
de $V$
tel que l'image de $\bar\Lambda$ dans $V/M$ ne contient pas de droite,
la proposition~\ref{prop:X-controle} implique donc
que la fonction $\mathsf X_\Lambda$
appartient à l'espace $\mathcal H_M(\Lambda;\Lambda)$
défini par la classe de contrôle $\mathscr D^{\max}$.

Le théorème principal de cette section est le suivant.
\begin{theo}\label{theo:prolongement}%
Soit $M\subset V$ un sous-espace vectoriel muni d'une mesure
de Lebesgue. 

Soit $C$ l'intersection de $\Lambda$ avec un voisinage
convexe de $0$ et soit $f\in\mathcal H_M(\Lambda;C)$.
Soit $M'$ un sous-espace vectoriel de $M$, $\pi$ la projection
$V\ra V'=V/M'$, $\Lambda'=\pi(\Lambda)$ et $C'=\pi(C)$.

Alors, la fonction $\mathscr S_{M'}(f)$
appartient à $\mathscr H_{M/M'}(\Lambda';C')$.

Si de plus l'adhérence du cône $\Lambda'$ ne contient pas de droite
et si pour tout $z\in\Lambda$,
\[ \lim_{s\ra 0^+ } \frac{f(sz)}{\mathsf X_\Lambda(sz)} = 1, \]
alors pour tout $z'\in\Lambda'$,
\[ \lim_{s\ra 0^+}  \frac{\mathscr S_{M'}(f) (sz')}{\mathsf X_{\Lambda'}(sz')}=
1. \]
\end{theo}

\begin{coro}
Supposons de plus que $f$ est la restriction à $\Lambda\cap C$
d'une fonction holomorphe $M$-contrôlée sur $\Lambda$.
Alors, la fonction $\mathscr S_M(f)$ sur $V'$ est méromorphe dans
un voisinage convexe de $\Lambda'$, ses pôles étant simples définis
par les faces (de codimension $1$) de $\Lambda'$.
\end{coro}

\subsection{D\'emonstration du th\'eor\`eme}

D'après le lemme~\ref{lemm:integrale}, la fonction $\mathscr S_{M'}(f)$
est holomorphe et $M/M'$-contrôlée sur $\Tube(C')$.
Le but est de montrer qu'elle y est la restriction d'une fonction
m\'eromorphe dont on contr\^ole les p\^oles et la croissance.
La démonstration est fondée sur l'application successive du
théorème des résidus pour obtenir le prolongement méromorphe.
La définition des classes de contrôle est faite pour assurer
l'intégrabilité ultérieure de chacun des termes obtenus.

Par récurrence, il suffit de démontrer le résultat lorsque
$\dim M'=1$. Soit $m_0$ un générateur de $M'$.
Munissons la droite $\R m_0$ de la mesure de Lebesgue $d\rho$.
Soit $\Phi\subset V^*$
une famille de formes linéaires deux à deux non proportionnelles
positives sur $\Lambda$ et 
dont les noyaux sont les faces de $\Lambda$.

Soit $B$ un ouvert convexe et symétrique par rapport à l'origine,
assez petit de sorte que pour tout
$\phi\in \Phi$ et tout $v\in B$, $\abs{\phi(v)}<1$ et que
la fonction
\[ g(z) = f(z) \prod_{\phi\in\Phi} \frac{\phi(z)}{1+\phi(z)} \]
admette un prolongement holomorphe $M$-contrôlé sur $\Tube(B)$.
L'intégrale à étudier est
\[ \int_{-\infty}^{+\infty} g(z+itm_0) \prod_{\phi\in\Phi}
\frac{1+\phi(z+itm_0)}{\phi(z+itm_0)}\, dt.  \]

On veut déplacer la droite d'intégration vers la gauche.
Fixons $\tau>0$ tel que $ 2\tau m_0 \in B$.
Ainsi, si $\Re(z)\in \frac12 B$,
$z+(u+it)m_0$
appartient à $\Tube(B)$
pour tout $u\in [-\tau;0]$ et tout $t\in\R$.

Notons $\Phi^+$, $\Phi^-$ et $\Phi^0$ les ensembles des $\phi\in\Phi$ tels que
respectivement $\phi(m_0)>0$, $\phi(m_0)<0$ et $\phi(m_0)=0$.
Soit $B_1\subset \frac12 B$ l'ensemble des $v\in \frac12 B$ tels que 
pour tout $\phi\in\Phi^+$,
$\abs{\phi(v)}<\frac\tau 2 \phi(m_0)$.

Dans la bande $-\tau \leq s \leq 0$,
les pôles de la fonction holomorphe
\[ s\mapsto g(z+sm_0)\prod_{\phi\in\Phi} \frac{1+\phi(z+sm_0)}{\phi(z+sm_0)} \]
sont ainsi donnés par
\[ s_\phi(z) = - \frac{\phi(z)}{\phi(m_0)}, \quad \phi\in\Phi^+. \]
Le pôle $s=s_\phi(z)$ est simple si et seulement si pour tout $\psi\in\Phi^+$
tel que $\psi\neq\phi$, 
\[ \phi(z) \psi(m_0) - \psi(z) \phi(m_0) \neq 0 . \]
Comme $\phi$ et $\psi$ sont non proportionnelles, 
$\psi(m_0)\phi-\phi(m_0)\psi$ est une forme linéaire non nulle ;
notons $B_1^\dagger\subset B_1$ le complémentaire des hyperplans qu'elles
définissent lorsque $\phi\neq\psi$ parcourent les éléments de $\Phi^+$.

Si $z\in\Tube(B_1^\dagger)$ et si $T> \max \{ \abs{\Im(s_\phi(z))}\,;\,
\phi\in\Phi^+\}$,
la formule des résidus pour le contour
délimité par le rectangle $-\tau\leq\Re(s)\leq 0$, $-T\leq\Im(s)\leq T$
s'écrit
\begin{align*}
\int_{-T}^T  g(z+itm_0) \prod_{\phi\in\Phi}
\frac{1+\phi(z+itm_0)}{\phi(z+itm_0)}\, dt \hskip -4cm \\
& = \sum_{\phi\in\Phi^+}
   \frac{2i\pi}{\phi(m_0)} g(z+s_\phi(z)m_0) \prod_{\psi\neq\phi}
        \frac{1+\psi(z+s_\phi(z)m_0)}{\psi(z+s_\phi(z)m_0)} \\
& \quad{} + \int_{-T}^T g(z-\tau m_0+itm_0) \prod_{\phi\in\Phi}
\frac{1+\phi(z-\tau m_0 +itm_0)}{\phi(z-\tau m_0 +itm_0)}\, dt  \\
& \quad{} + \int_{0}^{-\tau} 
       g(z+sm_0+iTm_0) \prod_{\phi\in\Phi}
 \frac{1+\phi(z+sm_0+iTm_0)}{\phi(z+sm_0+iTm_0)}\, ds \\
& \quad{} + \int_{-\tau}^{0} 
       g(z+sm_0-iTm_0) \prod_{\phi\in\Phi}
 \frac{1+\phi(z+sm_0-iTm_0)}{\phi(z+sm_0-iTm_0)}\, ds .
\end{align*}

Lorsque $T\ra +\infty$, l'hypothèse que $g$ est $M$-contrôlée
et l'axiome~(\ref{defi:classe},\ref{axiome:limite})
des classes de contrôles impliquent que
ces deux dernières intégrales (sur les segments horizontaux
du rectangle) tendent vers~$0$.
De même, l'axiome~(\ref{defi:classe},\ref{axiome:integrale})
assure la convergence des deux premières intégrales
vers les intégrales correspondantes de $-\infty $ à $+\infty$.

Par suite, si $z\in\Tube(B_1^\dagger\cap\Lambda)$, on a
\begin{multline}\label{eq:residus}%
\mathcal S_{\R m_0}(f)(z)
 = \sum_{\phi\in\Phi^+}
   g(z+s_\phi(z)m_0) \prod_{\psi\neq\phi}
        \frac{1+\psi(z+s_\phi(z)m_0)}{\psi(z+s_\phi(z)m_0)} \\
 + \prod_{\phi\in\Phi^0} \frac{1+\phi(z)}{\phi(z)}
    \, 
   \int_{-\infty}^\infty g(z-\tau m_0+itm_0)
  \prod_{\phi\in\Phi\setminus\Phi^0}
  \frac{1+\phi(z-\tau m_0 +itm_0)}{\phi(z-\tau m_0 +itm_0)}\, dt.
\end{multline}
Il résulte alors des axiomes~(\ref{defi:classe},\ref{axiome:projecteur})
et~(\ref{defi:classe},\ref{axiome:integrale})
des classes de contrôles que la fonction
\begin{equation}\label{eq:poles1}%
z\mapsto \mathcal S_{\R m_0}(f) (z)
\prod_{\phi\in\Phi^0} \frac{\phi(z) }{1+\phi(z)}
\prod_{\phi\in\Phi^+} \prod_{\psi\not\in\Phi^0\cup\{\phi\}} 
       \frac{\phi(s+s_\psi(z)m_0)}{1+\psi(s+s_\phi(z)m_0)}
\end{equation}
définie sur $\Tube(B_1^\dagger\cap\Lambda)$ s'étend en une
fonction holomorphe $M/M'$-contrôlée sur $  \Tube(\pi(B_1^\dagger))$.
En particulier, $\mathcal S_{\R m_0}(f) $ se prolonge méromorphiquement
à $\Tube(B_1^\dagger)$
et les pôles de $\mathcal S_{\R m_0} (f)$
sont donnés par une famille finie de formes linéaires.
Le lemme suivant les interprète géométriquement.

\begin{lemm}
Les faces de $\Lambda'$ sont les noyaux des formes linéaires
deux à deux non proportionnelles
sur $V/\R m_0$
$\phi\in\Phi^0$ et $\phi-\frac{\phi(m_0)}{\psi(m_0)}\psi$
pour $\phi\in\Phi^+$ et $\psi\in\Phi^-$.

De plus, si $\phi$ et $\psi\in\Phi^+$, le noyau de
$\phi-\frac{\phi(m_0)}{\psi(m_0)}\psi$ rencontre $\Lambda'$.
\end{lemm}
\begin{proof}
Un vecteur $x\in V$ appartient à $\Lambda$ si et seulement
si $\phi(x)>0$ pour tout $\phi\in\Phi$.
Par suite, $\pi(x)\in\Lambda'$ si et seulement si il
existe $\alpha\in\R$ tel que $\phi(x-\alpha m_0)>0$ pour tout
$\phi\in\Phi$.
Si $\phi\in\Phi^0$, cette condition est exactement $\phi(x)>0$.
Pour les autres $\phi$, elle devient
\[ \max_{\phi\in\Phi^-} \frac{\phi(x)}{\phi(m_0)}
< \alpha < \min_{\phi\in\Phi^+} \frac{\phi(x)}{\phi(m_0)} \]
d'où la première partie du lemme.

Pour la seconde, soit $\phi$ et $\psi$ deux éléments distincts
de $\Phi^+$. Si le noyau de $\phi-\frac{\phi(m_0)}{\psi(m_0)}\psi$
ne recontre pas $\Lambda'$, quitte à permuter $\phi$ et $\psi$,
on a 
\[ \frac{\phi(v)}{\phi(m_0)} > \frac{\psi(v)}{\psi(m_0)}  \]
pour tout $v\in\Lambda$ et cela contredit le fait que $\phi$
et $\psi$ définissent deux faces distinctes de $\Lambda$.
\end{proof}

On sait que $\mathscr S_{\R m_0} (f)$
est holomorphe sur $\Tube(\Lambda')$.
Il résulte du lemme que les formes linéaires 
$ \psi+s_\phi(z)\phi$ avec $\phi\in\Phi^+$ et
$\psi\not\in\Phi^0\cup\{\phi\}$ 
sont des pôles apparents dès que $\psi\in\Phi^+$.
Les autres correspondent aux faces de $\Lambda'$ !

Autrement dit, nous avons déjà prouvé que $\mathscr S_{\R m_0}(f)$
est la restriction à $\Tube(\pi(B_1))$ d'une fonction méromorphe dont les pôles
(simples) sont donnés par les faces de $\Lambda'$.
Montrons comment contrôler la croissance de $\mathscr S_{\R m_0}(f)$ dans les
bandes verticales.

\begin{lemm}
Soit $V$ un espace vectoriel, $M$ un sous-espace vectoriel,
$B$ un voisinage de $0$ dans $V$.
Soit $h$ une fonction holomorphe sur $\Tube(B)$
et soit $\ell$ une forme linéaire sur $V$.
Si la fonction $z\mapsto h(z) \frac{\ell(z)}{1+\ell(z)}$ est $M$-contrôlée,
$h$ est $M$-contrôlée.
\end{lemm}
\begin{proof}
Il faut montrer que $h$ est $M$-contrôlée dans un voisinage de tout
point de $B$.
Soit donc $x_0\in B$ et $K$ un voisinage compact de $x_0$ contenu dans $B$.
Soit $\kappa\in\mathscr D(M,V)$ telle que pour tout $x\in K$
et tout $y\in V$,
\[ \abs{ h(x+iy) \frac{\ell(x+iy)}{1+\ell(x+iy)}} \leq \kappa(y). \]

Supposons d'abord que $\ell(x_0)\neq 0$.
Si $\rho=\abs{\ell(x_0)}/2>0$,
il existe un voisiange compact $K_1\subset K$ de $x_0$
où $\abs{\ell}\geq\rho$.
Alors, pour tout $x\in K_1$ et tout $y\in V$, on a 
\[ \abs{h(x+iy)} \leq \kappa(y) \frac{1+\abs{\ell(x+iy)}}{\ell(x+iy)}
 \leq \frac{1+\rho}{\rho}\kappa(y) , \]
ce qui prouve que  $h$ est $M$-contrôlée dans $K_1$.

Si $\ell(x_0)=0$, soit $u\in V$ tel que $\ell(u)=1$, $K_1$ un voisinage
compact de $x_0$ assez petit et $\rho>0$
tels que pour tout $t\in [-1;1]$ et tout $x\in K_1$,  $x+t\rho u\in K$.
La fonction $s\mapsto h(x+iy+s\rho u)$ est une fonction holomorphe 
sur le disque unité fermé $\abs{s}\leq 1$.
D'après le principe du maximum, on a donc pour tout $x+iy\in \Tube(K_1)$,
\[ \abs{h(x+iy)}\leq \sup_{\abs{s}\leq 1} 
\abs{h(x+iy+s\rho u)}
 = \sup_{\abs{s}=1} \abs{h(x+iy+s\rho u)}
 \leq \frac{1+\rho}\rho \sup_{\abs{s}\leq 1} \kappa(y+su) . \]
L'axiome~(\ref{defi:classe},\ref{axiome:max}) assure alors l'existence d'une
fonction $\kappa_1\in \mathcal D(M,V)$
telle que pour tout $x+iy\in\Tube(K_1)$,
\[ \abs{h(x+iy)} \leq \kappa_1(y). \]
La fonction $h$ est donc $M$-contrôlée dans un voisinage de $x_0$.
\end{proof}

Il reste à démontrer que si pour tout $z\in \Lambda$,
$\lim\limits_{t\ra 0^+} f(tz) /  \mathsf X_\Lambda(tz) = 1$,
alors 
\[ \lim\limits_{t\ra 0^+} \mathscr S_{\R m_0}(f)(tz') / \mathsf
X_{\Lambda'}(tz')=1 . \]

Comme $X_\Lambda(tz)=t^{-\dim V}X_\Lambda(z)$,
l'hypothèse $f(tz)\sim \mathsf X_\Lambda(tz)$ se récrit
\[ \lim_{t\ra 0} t^{\dim V-\#\Phi} g(tz) = \mathsf X_\Lambda(z). \]
D'autre part, la formule~\eqref{eq:residus} donne
\begin{align*}
t^{-1+\dim V} \mathscr S_{\R m_0}(f) (tz) \hskip -2cm \\ 
&  = t^{-1+\dim V} \sum_{\phi\in\Phi^+}
   g(tz+s_\phi(tz)m_0) \prod_{\psi\neq\phi}
        \frac{1+\psi(tz+s_\phi(tz)m_0)}{\psi(tz+s_\phi(tz)m_0)} \\
& \qquad{} + t^{-1+\dim V} \prod_{\phi\in\Phi^0} \frac{1+\phi(tz)}{\phi(tz)} 
    \times  \\
   & \qquad\qquad {}\times
   \int_{-\infty}^\infty g(tz-\tau m_0+itm_0)
  \prod_{\phi\in\Phi\setminus\Phi^0}
  \frac{1+\phi(tz-\tau m_0 +itm_0)}{\phi(tz-\tau m_0 +itm_0)}\, dt \\
& = \sum_{\phi\in\Phi^+} t^{\dim V-\#\Phi} g(t(z+s_\phi(z)m_0))
        \prod _{\psi\neq\phi}
        \frac{1+t\psi(z+s_\phi(z)m_0)}{\psi(z+s_\phi(z)m_0)} \\
& \qquad{} + t^{-1+\dim V-\#\Phi^0}
   \prod_{\phi\in\Phi^0} \frac{1+t\phi(tz)}{\phi(z)} \times \\
   & \qquad\qquad {}\times
   \int_{-\infty}^\infty g(tz-\tau m_0+itm_0)
  \prod_{\phi\in\Phi\setminus\Phi^0}
  \frac{1+\phi(tz-\tau m_0 +itm_0)}{\phi(tz-\tau m_0 +itm_0)}\, dt .
\end{align*}
Un vecteur non nul de $V$ ne peut appartenir qu'à au plus $\dim V-1$
faces de $\Lambda$ et seuls les générateurs de $\Lambda$
appartiennent à $\dim V-1$ faces.
Comme $m_0$ est supposé n'être pas un générateur de $\Lambda$,
$\#\Phi^0\leq \dim V-2$.
Lorsque $t$ tend vers~$0$, on a donc
\[ \lim t^{-1+\dim V}\mathscr S_{\R m_0}(f) (tz)
= \sum_{\phi\in\Phi^+} \mathsf X_\Lambda(z+s_\phi(z)m_0)
  \prod_{\psi\neq\phi} \frac{1}{\psi(z+s_\phi(z)m_0)}  \]
où le second membre ne dépend plus de $f$. Comme on peut appliquer cette
formule à $f=\mathsf X_\Lambda$, on obtient donc
\begin{align*}
 \lim t^{1-\dim V}(\mathscr S_{\R m_0} (f))(tz)
     = \lim t^{1-\dim V}(\mathscr S_{\R m_0}(\mathsf X_\Lambda)) (tz) \\
     = \lim t^{1-\dim V}\mathsf X_{\Lambda'}(tz)
       = \mathsf X_{\Lambda'}(z). 
\end{align*}
Le théorème est ainsi démontré.

\begin{rema}
La d\'emonstration s'adapte sans peine lorsque $f$ d\'epend uniform\'ement
de param\`etres suppl\'ementaires.
\end{rema}


\section{Vari\'et\'es toriques}
\label{sec:torique}%

\def\PL{\mathrm{PL}}

Dans ce paragraphe, nous 
montrons comment les raffinements
analytiques du paragraphe~\ref{sec:analyse} permettent de préciser
le développement asymptotique
obtenu par Batyrev--Tschinkel dans~\cite{batyrev-t98b}
pour la fonction zêta des hauteurs d'une variété torique.
Les résultats techniques que nous rappelons à l'occasion
seront réutilisés au paragraphe suivant, lorsque nous
traiterons le cas d'une fibration en variétés toriques.

\subsection {Pr\'eliminaires}

\paragraph{Rappels adéliques}
\label{rappel:caracteres}%
Notons $S=\Spec\mathfrak o_F$
le spectre de l'anneau des entiers de $F$.
Si $v$ est une place de $F$, on définit la norme $\norm{\cdot}_v$
sur $F_v$ de la manière habituelle,
comme le \emph{module} associé à une mesure de Haar
additive sur $F_v$. En particulier, si $\pi_v$ est une uniformisante
en une place finie $v$, $\norm{\pi_v}_v$ est l'inverse
du cardinal du corps résiduel en~$v$.

Soit $G$ un tore d\'eploy\'e de dimension~$d$ sur $S$.
D\'esignons par $\mathbf K_\infty$ la collection
de ses sous-groupes compacts maximaux aux places \`a l'infini
et $\mathbf K_G=\prod_{v\nmid\infty}G(\mathfrak
o_v)\prod_{v\mid\infty}\mathbf K_v\subset G(\mathbf A_F)$.
Il nous faut faire quelques rappels sur
la structure du groupe $\mathcal A_G$ des caract\`eres 
de $G(F)\backslash G(\A_F)/\mathbf K_G$.
On a un homomorphisme de noyau fini
$\mathcal A_G\ra \bigoplus_{v|\infty} X^*(G)_\R$,
$\chi\mapsto\chi_\infty$
obtenu en associant \`a un caract\`ere ad\'elique son \emph{type \`a l'infini},
c'est-\`a-dire sa restriction
au sous-groupe de $G(\A)$ dont les composantes aux places
finies sont triviales.
En choisissant une norme sur $X^*(G)_\R$, on obtient 
ainsi une {\og norme\fg} $\chi\mapsto \norm{\chi_\infty}$
sur $\mathcal A_G$.

Il existe enfin un homomorphisme $X^*(G)_\R\ra \mathcal A_G$,
tel que l'image du caract\`ere alg\'ebrique $\chi\in X^*(G)$
est le caract\`ere ad\'elique $\mathbf g\mapsto \abs{\chi(\mathbf g)}^i$
dont le type \`a l'infini s'identifie \`a $\chi$ sur chaque
composante.

Le quotient $\mathcal A_G/X^*(G)_\R$ est un $\Z$-module de type
fini et de rang $(\rho-1)d$ (o\`u $\rho=r_1+r_2$, $r_1$ et $r_2$
d\'esignant comme d'habitude les nombres de places r\'eelles et complexes)
et l'on peut fixer une d\'ecomposition
$\mathcal A_G = X^*(G)_\R \oplus \mathcal U_G$, par
exemple \`a l'aide d'un scindage de la suite exacte
\[ 1 \ra \gm(\A_F)^1 \ra \gm(\A_F) \xrightarrow{\abs{\cdot}} \R^* \ra 1. \]
(Rappelons que $G$ est suppos\'e d\'eploy\'e.)

\paragraph{Rappels sur les vari\'et\'es toriques}

Notons $M=X^*(G)_\R$, c'est un espace vectoriel sur $\R$
de dimension finie~$d$.
Consid\'erons une compactification \'equivariante $\mathcal X$
de $G$, lisse sur $S$.
D'apr\`es la th\'eorie des vari\'et\'es toriques (cf.\ par exemple~\cite{oda88},
\cite{fulton93}), $\mathcal X$ est d\'efinie par un \'eventail 
complet et régulier~$\Sigma$
de~$N:=\Hom(M,\R)$ form\'e de c\^ones convexes simpliciaux rationnels.
Il existe ainsi une famille (minimale)
$(e_j)_{j\in J}$ de vecteurs de $N$
telle que tout c\^one $\sigma\in\Sigma$ soit engendr\'e par
une sous-famille $(e_j)_{j\in J_\sigma}$ de cardinal
$\dim\operatorname{vect}(\sigma)$.
On note $\Sigma(d)$ l'ensemble des c\^ones de $\Sigma$
de dimension $d$.

L'espace vectoriel $\mathrm{PL}(\Sigma)$
des fonctions continues $N\ra\R$ dont 
la restriction \`a chaque c\^one de $\Sigma$ est lin\'eaire
est un espace vectoriel de dimension finie sur $\R$,
d'ailleurs \'egale \`a $\# J$ ;
munissons le d'une norme arbitraire.
L'espace vectoriel $\Pic^G(\mathcal X_F)_\R$ est 
isomorphe \`a $\mathrm{PL}(\Sigma)$ ; il poss\`ede une base
canonique form\'ee des fibr\'es en droites $G$-lin\'earis\'es associ\'es aux
diviseurs $G$-invariants sur $\mathcal X_F$. \`A chaque
$e_j$ correspond un tel diviseur $D_j$ ; \`a un diviseur $G$-invariant
$D=\sum_{j} \lambda_j D_j$ correspond l'unique fonction
$\phi\in\mathrm{PL}(\Sigma)$ telle que $\phi(e_j)=\lambda_j$.
Dans cette description,
le c\^one des diviseurs effectifs correspond simplement l'ensemble
des \'el\'ements de $\Pic^G(\mathcal X_F)$ dont les coordonn\'ees
$(\lambda_j)$ v\'erifient $\lambda_j\geq 0$ pour tout $j$.
Plus généralement, on notera $\Lambda_t$ l'ensemble des éléments
de $\Pic^G(\mathcal X_F)$ tels que $\lambda_j>t$ pour tout $j$ ;
le cône ouvert $\Lambda_0$ est aussi noté $\mathrm{PL}^+(\Sigma)$
et encore $\Lambda_\eff^0(\mathscr X_F)$.

Cette base $(D_j)$ de $\Pic^G(\mathcal X_F)$ et l'homomorphisme canonique
$\iota:X^*(G)\ra \Pic^G(\mathcal X)$
induisent des sous-groupes \`a un param\`etre
$\gm\ra  G$, d'o\`u, pour tout caract\`ere $\chi\in\mathcal A_G$,
des caract\`eres $\chi_j$ de $\gm(F)\backslash \gm(\A_F)/\mathbf K_{\gm}$,
autrement dit des caract\`eres de Hecke.

Les fibr\'es en droites sur $\mathcal X_F$ seront syst\'ematiquement
munis de leur m\'etrique ad\'elique canonique introduite notamment
dans~\cite{batyrev-t95b}.
Cela nous fournit un homomorphisme canonique
$\Pic(\mathcal X_F)\ra\hatPic(\mathcal X)$
qui induit un homomorphisme 
\begin{equation}\label{eq:metcan}%
\Pic^G(\mathcal X_F)\ra \hatPic^{G,\mathbf K}(\mathcal X). 
\end{equation}
On v\'erifie ais\'ement, par exemple sur les formules donn\'ees
dans~\cite{batyrev-t95b}, que les sous-groupes compacts maximaux aux places
archim\'ediennes agissent de mani\`ere isom\'etrique.
De plus,  le choix d'une $G$-lin\'earisation
fournit une unique $F$-droite de sections ne s'annulant pas sur $G$,
donc en particulier une fonction hauteur sur les points ad\'eliques de
$\mathcal X_F$ comme dans la d\'efinition~\ref{hauteuradelique}.
Cette fonction s'\'etend en une application {\og bilin\'eaire\fg}
\[ H : \mathrm{PL}(\Sigma)_\C \times G(\A_F) \ra \C^* . \]
(On a identifi\'e $\Pic^G(\mathcal X_F)_\C $ et $\mathrm{PL}(\Sigma)_\C$.)

\begin{lemm}\label{lemm:chi/H}
Soit $m\in X^*(G)$ et notons $\chi_m\in\mathcal A_G$ le caract\`ere
ad\'elique qu'il d\'efinit. On a alors
\[ \chi_m(\mathbf g)  = H (\iota(m),\mathbf g)^{-i}. \]
\end{lemm}
\begin{proof}
Par d\'efinition, $\iota(m)$ est le fibr\'e en droite trivial
sur $\mathcal X$ muni de la $G$-lin\'earisation dans laquelle $G$
agit par multiplication par le caract\`ere alg\'ebrique $m$.
Ainsi, la droite de sections rationnelles $G$-invariante
et ne s'annulant pas sur $G$ est engendr\'ee par le caract\`ere $m$
vu comme fonction rationnelle sur $\mathcal X$.
La d\'efinition de $H$ implique que
\[ H(\iota(m),\mathbf g) = \prod_v \norm{m(g_v)}^{-1}
= \norm{m(\mathbf g)}^{-1}. \]
Or, 
\[ \chi_m( \mathbf g) = \norm{m(\mathbf g)}^i = H(\iota(m),\mathbf g)^{-i}.\]
\end{proof}

\paragraph{Mesures}
Pour toute place $v$ de $F$, on fixe une mesure de Haar $dx_v$ sur $F_v$.
On suppose que pour presque toute place finie $v$, la mesure 
du sous-groupe compact $\mathfrak o_v$ est égale à~$1$.
Alors, $dx=\prod_v dx_v$ est une mesure de Haar sur le groupe
localement compact $\A_F$.
On en déduit pour tout $v$ une mesure de Haar $\mu'_{\gm,v}=
\norm{x_v}_v^{-1} dx_v$ sur $F_v^*$.
Pour presque toute place finie~$v$, la mesure de $\mathfrak o_v^*$ est égale
à $1-q_v^{-1}$ ; définissons ainsi, si $v$ est une place finie,
$\mu_{\gm,v}=(1-q_v^{-1})^{-1}\mu'_{\gm,v}$.
On munit alors $\A_F^*$ de la mesure
\[  \prod_v \mu_{\gm,v} = 
   \prod_{v\nmid\infty}  (1-q_v^{-1})^{-1} \norm{x_v}^{-1} dx_v  \times
   \prod_{v|\infty} \norm{x_v}^{-1}dx_v . \]
Remarquons que  $\zeta_{F,v}(1)=(1-q_v^{-1})^{-1}$ est le facteur
local en la place finie $v$ de la fonction zêta de Dedekind du corps~$F$.

Tout $\mathfrak o_F$-isomorphisme $G\simeq \gm^d$ 
induit alors des mesures
de Haar $\mu'_{G,v}$ et $\mu_{G,v}=\zeta_{F,v}(1)^d \mu'_{G,v}$
sur $G(F_v)$ pour toute place $v$ de $F$, indépendantes
de l'isomorphisme. On en déduit aussi une mesure de Haar $\prod \mu_{G,v}$
sur $G(\A_F)$.

D'autre part, le fibré canonique sur $\mathscr X$ est métrisé.
Peyre a montré dans~\cite{peyre95} comment en déduire
une mesure sur $\mathscr X(\A_F)$.
Pour toute place $v$, on dispose d'une mesure
$\mu'_{\mathscr X,v}$ sur $\mathscr X(F_v)$ définie par la formule
\[ \mu'_{\mathscr X,v}  = \norm{d\xi_1\wedge \dots\wedge d\xi_d}_v^{-1}
          d\xi_1\,\dots\, d\xi_d \]
si $(\xi_1,\dots,\xi_d)$ est un système arbitraire de coordonnées locales
sur $\mathscr X(F_v)$.
Si l'on restreint la mesure $\mu'_{\mathscr X,v}$ à $G(F_v)$,
on obtient  donc
\begin{equation}\label{prop:mestam}
  H_v(-\rho,x) \mu'_{G,v} , 
\end{equation}
$\rho$ désignant la fonction de $\PL(\Sigma)$ telle que 
pour tout $j$, $e_j\mapsto 1$
($\rho$ correspond à la classe anticanonique).

Pour presque toute place finie $v$, on a alors
\[ \mu'_{\mathscr X,v}(\mathscr X(F_v)) = q_v^{-d} \# \mathscr X(k_v) . \]
La décomposition cellulaire des variétés toriques
(point n'est besoin ici d'invoquer le théorème de Deligne sur
les conjectures de Weil)
implique alors que 
\[ \# \mathscr X(k_v) = q_v^{d} + \rg(\Pic\mathscr X_F)
          q_v^{d-1} + O(q_v^{d-2}) . \]
Par suite, le produit infini
\[ \prod_{v\nmid\infty}
     \mu'_v (\mathscr X(F_v)) \zeta_{F,v}(1)^{-\rg (\Pic\mathscr X_F)} \]
est convergent.
Définissons une mesure $\mu_{\mathscr X,v}$ sur $\mathscr X(F_v)$ par
\[ \mu_{\mathscr X,v}= 
 \zeta_{F,v}(1)^{-\rg\Pic \mathscr X_F}\mu'_{\mathscr X,v}  \]
si $v$ est finie et $\mu_{\mathscr X,v}=\mu'_{\mathscr X,v}$
si $v$ est archimédienne.
Ainsi,
le produit infini
$\prod_v \mu_{\mathscr X,v}$ converge et définit
une mesure, dite \emph{mesure de Tamagawa} sur $\mathscr X(\A_F)$.
Le nombre de Tamagawa de $\mathscr X(\A_F)$ est alors définie par
\begin{equation}\label{defi:mestam}
 \tau(\mathscr X) =
    \mu(\A_F/F)^{-d}  \res_{s=1}\zeta_F(s)^{\rg(\Pic\mathscr X_F)} 
 \mu_{\mathscr X}(\A_F). 
\end{equation}

\begin{rema}
La différence de formulation avec la définition que donne
Peyre dans~\cite{peyre95} n'est qu'apparente. Peyre a choisi
la mesure sur $F_v$ de la façon suivante : si $v$ est une place finie,
$dx_v(\mathfrak o_v)=1$, si $v$ est une place réelle, $dx_v$
est la mesure de Lebesgue usuelle sur $\R$ et si $v$ est une place
complexe, $dx_v$ est le double de la mesure usuelle sur $\C$.
Le volume de $\A_F/F$ est alors égal à $\Delta_F^{1/2}$.
\end{rema}

\subsection{Transformations de Fourier}

On s'int\'eresse \`a la transform\'ee de Fourier de la fonction
$g\mapsto H(-\lambda,g)$ sur le groupe ab\'elien localement compact $G(\A_F)$.
Rappelons qu'on a noté $\Lambda_1$ l'ensemble des $\lambda\in\PL(\Sigma)$
tels que $\lambda_j>1$ pour tout $j$. 
Alors, si $\lambda\in\Tube(\Lambda_1)$,
la fonction $g\mapsto H(-\lambda,g)$ est int\'egrable
(cf.~\cite{strauch-t99}, \S\,3.4), si bien
que la transform\'ee de Fourier existe pour tout
$\lambda\in\Tube(\Lambda_1)$.
Elle se d\'ecompose par construction en un produit
$\Fourier H= \Fourier H_f\times\Fourier H_\infty$, o\`u
\[
\Fourier H_f = (\res_{s=1}\zeta_F(s) ) ^{-d} \prod_{v\nmid\infty}
         (1-q_v^{-1})^{-d}\Fourier H_v
\]
et
$\Fourier H_\infty= \prod_{v|\infty}
\Fourier H_v $
sont les produits des int\'egrales
locales (renormalisées) aux places finies et archim\'ediennes.
(Les transformées de Fourier locales existent même dès que
pour tout $j$, $\Re(\lambda_j)>0$.)

\begin{lemm}\label{lemm:hatHfini}%
Soit $\Lambda_{2/3}\subset\PL(\Sigma)$ la partie convexe
d\'efinie par $\lambda_j>2/3$ pour tout $j$.
Il existe une fonction 
\[c_f:\Tube(\Lambda_{2/3})\times \mathcal A_G \ra\C, 
  \quad (\lambda,\chi)\mapsto c_f(\lambda,\chi), \]
holomorphe en $\lambda$ telle que $\log\abs{c_f}$ est born\'ee
et telle que le produit des transform\'ees
de Fourier locales aux places non archim\'ediennes
s'\'ecrive, pour tout $\chi\in\mathcal A_G$ et tout $\lambda\in\Tube(\Lambda)$
\[ 
\Fourier H_f (-\lambda,\chi) = c_f(\lambda,\chi) \prod_{j} L(\lambda_j,\chi_j). 
\]
\end{lemm}
\begin{proof}
Si $\chi$ est fixé, c'est la proposition 2.2.6 de~\cite{batyrev-t95b}.
Le fait que $\log\abs{c_f}$ soit borné indépendamment de $\chi$
se déduit immédiatement de la preuve dans \emph{loc. cit.}
\end{proof}

\begin{coro}
La fonction $\Fourier H_f$ se prolonge en une fonction
m\'eromorphe pour $\lambda\in\Tube(\Lambda_{2/3})$.
Plus pr\'ecis\'ement,
le produit $\prod_j (\lambda_j-1)\Fourier H_f (-\lambda,\chi)$
se prolonge en une fonction holomorphe dans $\Tube(\Lambda_0)$
et
\[ \lim_{\lambda\ra (1,\ldots,1)} \prod_j (\lambda_j-1)\Fourier H_f
(-\lambda,\chi) = 0 \]
si et seulement si $\chi\neq\mathbf 1$.
\end{coro}

Comme cons\'equence facile de l'estimation par Rademacher des valeurs
des fonctions $L$ de Hecke pour les caract\`eres non ramifi\'es,
estimation qui repose sur le principe de Phragm\'en--Lindel\"of,
on obtient la majoration suivante :
\begin{coro}\label{coro:rademacher}%
Pour tout $\eps>0$,  
il existe $0<\delta<1/3$ et un réel $c_\eps$
tels que si $\Re(\lambda_j)> 1-\delta$,
\[
 \prod_{j} \frac{\abs{\lambda_j-1}}{\abs{\lambda_j}}
 \Fourier H_f(-\lambda,\chi) 
 \leq c_\eps
         \big( 1+ \norm{\Im(\lambda)} \big)^{\eps}
         \big( 1+ \norm{\chi}\big)^{\eps}.
\]
\end{coro}

Passons maintenant aux places archim\'ediennes.
La proposition suivante précise 
la proposition 2.3.2 de~\cite{batyrev-t95b}.
\begin{prop}\label{prop:fourierinfini}%
Pour tout compact
$K\subset \Lambda_{2/3}\subset\mathrm{PL}(\Sigma)_\R$,
il existe un réel
$c_K$ telle que pour tout $\phi\in \Tube(K)$ et tout $m\in M$,
on ait la majoration 
\[ \abs{\mathscr F(m)} \leq c_K
     \frac{1}{1+\norm{m}} \sum_{\sigma\in\Sigma(d)}
           \frac{1+\norm{\phi}_\sigma}{
  \prod_{j\in J_\sigma} \left( 1+ \abs{\langle e_j,m\rangle}\right)}. \]
\end{prop}

\begin{coro}\label{coro:fourierinfini}%
Désignons par $\tilde\Sigma$ l'éventail $\prod_{v|\infty} \Sigma$
dans $\tilde N = \prod_{v|\infty} N$.
Si $\phi\in\PL(\Sigma)$, désignons par $\tilde\phi$
la fonction $\tilde N\ra\R$ définie par $(n_v)_v\mapsto \sum \phi(n_v)$.
Pour tout compact $K$ de $\PL(\Sigma)$ contenu dans $\Lambda_{2/3}$,
il existe une constante $c_K$ telle que pour tout $\phi\in\Tube(K)$
et tout $\chi\in\mathcal A_G$ décomposé sous la forme
$\chi=im+\chi_u\in iM\oplus\mathcal U_G$, on ait
\[ \check H_\infty(\phi,\chi) \leq \frac{c_K}{1+\norm{\chi}}
     \sum_{\tilde\sigma\in\tilde\Sigma}
     \frac{1+\norm{\Im\tilde\phi}_{\tilde\sigma}}
        { \prod_{e\in\tilde\sigma} \left(
          1+\abs{\langle e, \Im \tilde\phi|_{\tilde\sigma}+\tilde m\rangle}
                                  \right)
        } . \]
\end{coro}
\begin{proof}
Si l'on note $\tilde m = (m_v)_v$ la décomposition de $\chi$ à l'infini,
on remarque que 
\[
\check H_\infty(\phi,\chi)  = \prod_{v|\infty} \check H_v(\phi,\chi)  
   = \prod_{v|\infty} \mathscr F(\phi,m_v) 
  = \mathscr F(\tilde\phi,\tilde m).
\]
Il suffit alors d'appliquer la proposition précédente.
\end{proof}

\begin{proof}[Preuve de la proposition~\ref{prop:fourierinfini}]
Il faut estimer
\[ \mathscr F(m) = \int_N \exp( -\phi(v)-i\langle v,m\rangle)\, dv .\]
Soit $\sigma\in\Sigma$ un cône de base $(e_1,\ldots,e_d)$.
Si $\abs{\det(e_j)}$ désigne la mesure du parralèlotope de base les $e_j$,
on a
\begin{align}
\int_\sigma \exp(-\phi(v)-i\langle v,m\rangle)\, dv
  &     = \int_{\R_+^d} \prod_{j=1}^d \exp\big(-t_j(\phi(e_j)+i \langle
e_j,m\rangle) \big)\,{\abs{\det(e_j)}}\, \prod dt_j  \notag \\
  &     = c(\sigma) \prod_{j=1}^d \frac{1}{\phi(e_j)+i\langle e_j,m\rangle}. 
\end{align}
Ainsi, on a
\begin{equation}\label{eq:finfini-0}%
 \mathscr F(m) = \sum_\sigma c(\sigma)
         \prod_{e\in\sigma} \frac{1}{\phi(e)+i\langle e,m\rangle}.
\end{equation}
D'autre part, supposons que $m_j\neq 0$, on peut intégrer par parties
et écrire
\begin{align}
\mathscr F(m) &= \int_N \frac{1}{im_j} \left( -\frac{\partial\phi}{\partial
v_j}\right) \exp(-\phi(v)-i\langle v,m\rangle)\, dv \notag \\
-i m_j \mathscr F(m) &= \int_N \left( \frac{\partial\phi}{\partial
v_j}\right) \exp(-\phi(v)-i\langle v,m\rangle)\, dv \notag \\
&= \sum_{\sigma} \left. \frac{\partial\phi}{\partial v_j}\right|_\sigma
            \int_\sigma \exp(-\phi(v)-i\langle v,m\rangle)\, dv \notag \\
&= \sum_\sigma c(\sigma)
              \left. \frac{\partial\phi}{\partial v_j}\right|_\sigma
              \prod_{e\in\sigma} \frac{1}{\phi(e)+i\langle e,m\rangle}.
\label{eq:finfini-j}%
\end{align}
En combinant les égalités~\eqref{eq:finfini-0} et~\eqref{eq:finfini-j}
pour tous les indices $j$ tels que $m_j\neq 0$, on obtient une majoration
\[
\abs{\mathscr F(m)} \leq \frac{1}{1+\norm{m}} \sum_\sigma c(\sigma)
           \frac{1+\norm{\phi}_\sigma}{\prod_{e\in\sigma}
                      \abs{\phi(e)+i\langle e,m\rangle}}. 
\]
Finalement, comme $\phi\in\Tube(K)$, 
on a une estimation
\[ \abs{\phi(e)+i\langle e,m\rangle}
\gg 1+\abs{\Im(\phi)(e)+\langle e,m\rangle} \]
et la proposition s'en déduit.
\end{proof}

\subsection{Définition d'une classe de contrôle}
\label{subsec:controle}%

Soit $\beta$ un réel strictement positif.
Si $M$ et $V$ sont deux $\R$-espaces vectoriels de dimension finie
avec $M\subset V$, notons $\mathscr D_{\beta,\eps}(M,V)$ le sous-monoïde
de $\mathscr F(V,\R_+)$ engendré 
par les fonctions $h:V\ra\R_+$ telles que pour tout $\eps>0$,
il existe $c>0$,
$\eps\in\left]0;1\right[$ et une famille $(\ell_j)$ de formes linéaires
sur $V$ vérifiant :
\begin{itemize}
\item la famille $(\ell_j|_M)$ forme une base de $M^*$ ;
\item pour tout $v\in V$ et tout $m\in M$, on a
\begin{equation}
h(v+m) \leq c \frac{(1+\norm v)^{\beta}}{(1+\norm{m})^{1-\eps}}
        \frac{1}{\prod(1+\abs{\ell_j(v+m)})}.
\end{equation}
\end{itemize}
Notons alors $\mathscr D_\beta=\bigcap_{\eps>0} \mathscr D_{\beta,\eps}$.

\begin{prop}\label{prop:controle}
Les $\mathscr D_\beta({M,V})$ définissent une classe de contrôle au
sens de la définition~\ref{defi:classe}.
\end{prop}

La preuve de cette proposition consiste en une série
d'inégalités faciles mais techniques. Nous la repoussons
à l'appendice~\ref{app:controle}.

\subsection{La fonction zêta des hauteurs et la formule de Poisson}

On s'intéresse fonction zêta des hauteurs de $\mathscr X$ 
restreinte à l'ouvert dense formé par le tore $G$ ;
c'est par définition la série génératrice
\[ Z(\lambda) = \sum_{x\in G(F)}  H(-\lambda,x) , \]
quand elle converge.
Des théorèmes taubériens standard (voir l'appendice)
permettront de déduire de résultats analytiques sur $Z$ un
développement asymptotique du nombre de points de hauteur bornée
\[ N(\lambda, H) =\# \{ x\in G(F)\,;\, H(\lambda,x)\leq H \}. \]

\begin{lemm}\label{lemm:convergence}
Lorsque $\Re(\lambda)$ décrit un compact de $\Lambda_1$,
la fonction zêta des hauteurs converge uniformémént en $\lambda$.
Plus généralement, la série
\[ \sum_{x\in G(F)} H(-\lambda,x\mathbf g) \]
converge absolument uniformément lorsque $\Re(\lambda)$ décrit
un compact de $\Lambda_1$ et $\mathbf g$ un compact de $G(\A_F)$.
\end{lemm}
\begin{proof}
Compte tenu d'estimations pour $H(-\lambda,x\mathbf g)/H(-\lambda,x)$
lorsque $\mathbf g$ décrit un compact de $G(\A_F)$,
$x\in G(F)$ et $\lambda\in\Tube(\Lambda_1)$,
c'est en fait un corollaire de l'intégrabilité de la fonction
$H(-\lambda,\cdot)$ sur $G(\A_F)$.
Voir~\cite{batyrev-t98b}, Th. 4.2 et aussi~\cite{strauch-t99}, Prop.~4.3.
\end{proof}

Par conséquent, on peut appliquer la formule sommatoire de Poisson
sur le tore adélique $G(\A_F)$ pour le sous-groupe discret
$G(F)$. Compte tenu  de l'invariance de l'accouplement
de hauteurs par le sous-groupe compact maximal $\mathbf K_G$
de $G(\A_F)$, on en déduit la formule
\begin{equation}\label{eq:poisson.a}
Z(\lambda) = \int_{\mathscr A_G} \check H(-\lambda,\chi)\, d\chi 
\end{equation}
où $d\chi$ est la mesure de Haar sur le
groupe $\mathscr A_G$ des caractères unitaires continus
sur le groupe $G(F)\backslash G(\A_F)/\mathbf K_G$
duale de la mesure de comptage sur $G(F)$.

Rappelons que l'on a décomposé le groupe $\mathscr A_G= M\oplus \mathscr
U_G$, où $\mathscr U_G$ est un groupe discret.
De plus, si $\chi=m\oplus \chi_u$, 
\[ \check H(-\lambda,\chi) = \check H(-\lambda-im,\chi_u) \]
si bien que
\[ Z(\lambda)
      = \int_{M} \left( \sum_{\chi_u\in\mathscr U_G} \check
H(-\lambda-im,\chi_u) \right)\, dm \]
où $dm$ est la mesure de Lebesgue sur $M$ telle que
$dm\,d\chi_u=d\chi$, $d\chi_u$ étant la mesure de comptage sur $\mathscr
U_G$.

\begin{lemm}\label{lemm:normalisations}
Si $d^0m$ est la mesure de Lebesgue sur $M$ définie par le
réseau $M$, on a 
\[ dm = \left( 2\pi \vol(\A_F/F)  \res_{s=1}\zeta_F(s)\right)^{-d} 
         \,  d^0m .   \]
\end{lemm}
\begin{proof}
Par multiplicativité, il suffit de traiter le cas $G=\gm$ et $d=1$.
Notons $\A^1_F$ le sous-groupe de $\A_F^*$ formé des $x$ tels que
$\norm x=1$.
La suite exacte
\[  1 \ra \A^1_F/F^* \ra \A_F^*/F^* \xrightarrow{\log\norm x} \R \ra 0 \]
permet de munir $\A^1_F/F^*$ de la mesure de Haar
$dx^1$ telle que $d^*x=dx^1\,d^0n$.
La suite exacte duale
\[ 1 \ra \R \ra (\A_F^*/F^*)^* \ra (\A^1_F/F^*)^* \ra 1 \]
et la discrétude du groupe des caractères de $\A^1_F/F^*$ permet
de munir $(\A_F^*/F^*)^*$ de le mesure $d^0m\, \sum$.
Avec ces normalisations, la constante devant la formule de Poisson
est
$(2\pi \vol(\A^1_F/F^*))^{-1}$.
Compte tenu des normalisations choisies,
le théorème classique selon lequel $\tau(\gm)=\tau(\ga)=1$,
cf. par exemple~\cite{weil82}, p.~116, devient
\[ \vol(\A^1_F/F^*) = \vol(\A_F/F) \res_{s=1} \zeta_F(s), \]
d'où le lemme.
\end{proof}

Soit $\rho=(1,\dots,1)\in\PL(\Sigma)$. On décale la fonction
zêta des hauteurs de $\rho$  : si $\lambda\in\PL(\Sigma)^+$,
\[ Z(\rho+\lambda) = 
   \int_{M}
 \left( \sum_{\chi_u\in\mathscr U_G} \check H(-\lambda-\rho-im,\chi_u) \right)
  \, dm \]
Soit $F$ la fonction $\PL(\Sigma)^+\ra\C$ définie par la série
\[ \lambda\mapsto  (\vol(\A_F/F) \res_{s=1}\zeta_F(s))^{-d}
         \sum_{\chi_u\in\mathscr U_G} \check H(-1-\lambda,\chi_u),  \]
de sorte que si $\lambda\in\PL(\Sigma)^+$,
\begin{equation}\label{eq:poisson.b}
   Z(\lambda+\rho) = \frac{1}{(2\pi)^d} \int_{M} F(\lambda+im)\, d^0m .
\end{equation}

\begin{prop}\label{prop:Fcontrole}
Si $\beta>1$, la fonction $F$ appartient à l'espace $\mathscr H_{M}
(\PL(\Sigma)^+)$ défini par la classe de contrôle
$\mathscr D_\beta$ du paragraphe~\ref{subsec:controle}.

De plus, pour tout $\lambda\in\PL(\Sigma)^+$,
\[ \lim_{s\ra 0} \frac{F(s\lambda)}{\mathsf
X_{\PL(\Sigma)^+}(s\lambda)}
     =  \tau(\mathscr X),
\]
le nombre de Tamagawa de~$\mathscr X$.
\end{prop}
\begin{proof}
On a vu que l'on pouvait écrire
\[ \check H(-\rho-\lambda,\chi)
   = c_f(\lambda+\rho,\chi) 
      \check H_\infty(-\rho\lambda,\chi)
  \prod_j L(\lambda_j+1,\chi_{j}). \]
Par suite, la fonction
\[ \lambda\mapsto \check H(-\rho-\lambda,\chi) \prod_j
\frac{\lambda_j}{\lambda_j+1} \]
admet un prolongement holomorphe pour $\Re(\lambda_j)>-1$.

De plus, il résulte des corollaires~\ref{coro:rademacher}
et~\ref{coro:fourierinfini}
que pour tout $\eps>0$, il existe $\delta<1/3$ tel que
si pour tout $j$ on a $\Re(\lambda_j)>-\delta$, alors
\[
 \abs{\check H(-\rho-\lambda,\chi) 
 \prod_{j} \frac{\lambda_j}{\lambda_j+1} }
     \ll 
      \frac{(1+\norm{\Im(\lambda)})^{1+\eps}}{(1+\norm{\chi_\infty})^{1-\eps}}
      \sum_{\tilde\sigma\in\tilde\Sigma(d)}
         \frac{1}
              {\prod_{e\in \tilde\sigma_1}
        (1+\abs{\langle e,\Im(\lambda)|_{\tilde\sigma}+\chi_\infty
\rangle})},
\]
formule dans laquelle
$\chi_\infty$ désigne l'image de $\chi$ par l'homomorphisme
de noyau fini {\og type à l'infini\fg}
$\mathscr A_G\ra M_\infty=\bigoplus_{v|\infty} M$.
Ainsi, on obtient un prolongement
holomorphe de la fonction
$\Phi : \lambda\mapsto F(\lambda) \prod_j \lambda_j/(1+\lambda_j)$
pour $\Re(\lambda_j)>-\delta$ si l'on prouve
que pour tout $\tilde\sigma\in\tilde\Sigma(d)$, la série
\[ \sum_{\chi_u\in\mathscr U_G}
       \frac{1}{(1+\norm{\chi_{u,\infty}})^{1-\eps}}
        \frac{1}{\prod_{e\in \tilde\sigma_1}
        (1+\abs{\langle e,\Im(\lambda)|_{\tilde\sigma}+\chi_{u,\infty}
\rangle})}
\]
converge localement uniformément en $\lambda$ si $\Re(\lambda_j)>-\delta$.
Fixons $\tilde\sigma\in\tilde\Sigma(d)$. Alors, 
lorsque $e\in\tilde\sigma_1$,
 les formes linéaires $\langle e,\cdot\rangle$ forment une
base de $M_\infty^*$.
Il est facile de remplacer la sommation sur le sous-groupe
discret $\mathscr U_{G,\infty}$ par une intégrale sur
l'espace vectoriel qu'il engendre, lequel est d'ailleurs un supplémentaire
de $M$ envoyé diagonalement dans $M_\infty$.
La convergence est alors une conséquence de la
proposition~\ref{prop:majoration}.

Pour obtenir l'assertion sur la croissance de $F$, 
il faut montrer que si $\beta>1$, $K$ est un compact de $\PL(\Sigma)^+$,
$\lambda\in\Tube(K)$ et $m\in M$, on a une majoration
\[ \abs{\Phi(\lambda+im)}
       \ll \frac{(1+\norm{\Im (\lambda)})^{\beta}}{(1+\norm{m})^{1-\eps}}
   \sum _{\alpha} \prod_{k} \frac{1}{1+\abs{\ell_{\alpha,k}(\Im(\lambda)+m)}} 
\]
où $\alpha$ parcourt un ensemble fini et où
pour tout $\alpha$, $\{ \ell_{\alpha,k}\}_k$ est
une base de $\PL(\Sigma)^*$.
Il nous
faut récrire un peu différemment la majoration de $\check H$ obtenue ci-dessus
en remarquant que si la forme des transformées de Fourier
aux places finies fournit le prolongement méromorphe,
la convergence de la série provient, elle, des estimations
archimédiennes.
On écrit ainsi
\[
 \check H(-\rho-\lambda-im,\chi_u) 
   = c_f(\rho+\lambda+im,\chi_u) \prod_j L(\lambda_j+1+im,\chi_{u,j})
          \check H_\infty (-\rho-\lambda,\chi_m\chi_u ) \]
et donc
\begin{multline*}
\abs{\check H(-\rho-\lambda-im,\chi_u)
     \prod_j \frac{\lambda_j+im}{1+\lambda_j+im} } \\
\ll
\frac{(1+\norm{\Im(\lambda)+m})^{\eps}(1+\norm{\chi_u})^{\eps} }
     {1+\norm{m+\chi_{u,\infty}}}
\sum_{\tilde\sigma\in\tilde\Sigma(d)}
     \frac{1+\norm{\Im(\lambda)}_{\tilde\sigma}}
       {\prod_{e\in\tilde\sigma_1}
        (1+\abs{\langle
e,\Im(\lambda)|_{\tilde\sigma}+m+\chi_\infty\rangle})}.
\end{multline*}
Par suite,
\[
\abs{\Phi(\lambda+im)}
\leq \sum_{\tilde\sigma\in\tilde\Sigma(d)}
   (1+\norm{\Im(\lambda)}_{\tilde\sigma})
    (1+\norm{\Im(\lambda)+m})^\eps G_{\tilde\sigma}(\Im(\lambda),m) \]
où $\Phi_{\tilde\sigma(\phi,m)}$ est défini par la série
\[ \Phi_{\tilde\sigma}(\phi,m)
   = \sum_{\chi_u\in\mathscr U_G}
    \frac{(1+\norm{\chi_u})^\eps}{1+\norm{m+\chi_{u,\infty}}}
         \prod_{e\in\tilde\sigma_1} \frac{1}{1+\abs{\langle
e,\tilde\phi|_{\tilde\sigma}+m+\chi_{u,\infty}\rangle}}. 
\]
On a la majoration
\[ 1+\norm{\chi_u} \leq 1+\norm{m+\chi_{u,\infty}}+\norm{m}
        \leq (1+\norm{m+\chi_{u,\infty}})(1+\norm{m}) \]
et comme précédemment, on remplace la sommation sur le sous-groupe
discret $\mathscr U_G$ par l'intégrale sur l'espace vectoriel qu'il
engendre.
La proposition~\ref{prop:majoration} fournit alors pour tout $\eps'>\eps$
une estimation 
\[ G_{\tilde\sigma}(\phi,m)
 \ll \frac{1}{(1+\norm{m})^{1-\eps'}} \sum_\alpha
       \prod_k
\frac{1}{1+\abs{\ell_{\alpha,k}(m+\phi|_{\tilde\sigma})}} \]
où $\{\ell_{\alpha,k}\}_k$ est une base de $M^*$
et $\phi|_{\tilde\sigma}$ l'élément de $M$ qui coïncide avec
$(\phi,\dots,\phi)\in\bigoplus_{v|\infty}\PL(\Sigma)$
sur le cône $\tilde\sigma$ de l'éventail $\tilde\Sigma$.
L'application $\phi\mapsto \ell_{\alpha,k}(\phi|_{\tilde\sigma})$
est une forme linéaire $\ell_{\tilde\sigma,\alpha,k}$
sur $\PL(\Sigma)$.
On a ainsi
\begin{align*}
\abs{G(\lambda+im)}  & \ll 
   \frac{     (1+\norm{\Im(\lambda)})(1+\norm{\Im(\lambda)+m})^\eps}
       { (1+\norm{m})^{1-\eps-\eps'}}
 \sum_{\tilde\sigma} \sum_\alpha
       \prod_k \frac{1}{1+\abs{\ell_{\tilde\sigma,\alpha,k}
(\Im(\lambda)+m)}} \\
&\ll 
     \frac{(1+\norm{\Im(\lambda)})^{1+\eps}}
     {(1+\norm{m})^{1-2\eps-\eps'}}
  \sum_{\tilde\sigma} \sum_\alpha
       \prod_k \frac{1}{1+\abs{\ell_{\tilde\sigma,\alpha,k}
(\Im(\lambda)+m)}} .
\end{align*}
Comme on peut prendre $\eps$ et $\eps'$ arbitrairement petits,
la contrôlabilité est établie.

Il reste à calculer la limite quand $s\ra 0$ par valeurs supérieures
de $F(s\lambda)/\mathsf X_{\PL(\Sigma)^+}(s\lambda)$.
Le cône $\PL(\Sigma)^+$ est simplicial et
\[  \mathsf X_{\PL(\Sigma)^+}(\lambda) = \frac{1}{\prod_j \lambda_j}. \]
Ainsi,
\[ \frac{F(\lambda)}{\mathsf X_{\PL(\Sigma)^+}(\lambda)} =
   (\vol(\A_F/F)\res_{s=1}\zeta_F(s))^{-d}    \prod_j \lambda_j
      \sum_{\chi_u\in\mathscr U_G} \check H(-\lambda-\rho,\chi_u).
\]
D'après ce qui précède, la série qui définit
$F$ converge uniformément pour $\Re(\lambda_j)>-\delta$ ;
cela permet de permuter sommation et limite, si bien que
\begin{multline*}
 \lim_{s\ra 0^+} \frac{F(s\lambda)}{\mathsf X_{\PL(\Sigma)^+}(s\lambda)}
\\ = (\vol(\A_F/F)\res_{s=1}\zeta_F(s))^{-d}
  \sum_{\chi_u\in\mathscr U_G} \left(
 \lim_{s\ra 0^+} \check H(-s\lambda-\rho,\chi_u) \prod_j (s\lambda_j)
\right). 
\end{multline*}
En écrivant,
\[  \check H(-s\lambda-\rho,\chi) \prod_j (s\lambda_j)
      = c_f(s\lambda,\chi)
       \prod_j \big(s\lambda_j L(s\lambda_j+1,\chi_j)\big)
           \check H_\infty(s\lambda,\chi), \]
on voit que 
la limite est nulle si l'un des $\chi_j\neq\mathbf 1$
(car une des fonctions $L(\cdot,\chi_j)$ n'a pas de pôle en $1$,
les autres ont au plus un pôle simple).
\'Etudions maintenant le cas $\chi=\mathbf 1$.
Utilisant la formule~\eqref{prop:mestam}, il vient
\begin{align*}
\check H(-s\lambda-\rho,\mathbf 1) 
 \prod_j \zeta_F(1+\lambda_j s)^{-1}  \hskip -110pt \\
& = \prod_{v\nmid\infty} \zeta_v(1)^{d} \prod_j \zeta_v(1+\lambda_j s)^{-1}
      \int_{G(F_v)} H(-s\lambda-\rho,x) \mu'_{G,v}   \\
& \qquad\qquad {}\times
  \prod_{v|\infty}   \int_{G(F_v)} H(-s\lambda-\rho,x) \mu'_{G,v} \\
&  = \prod_{v\nmid\infty}
   \zeta_v(1)^d \prod_j \zeta_v(1+\lambda_j s)^{-1}
   \int_{G(F_v)} H(-s\lambda)\mu'_{\mathscr X,v}  \times
   \prod_{v|\infty} \int_{G(F_v)} H(-s\lambda)\mu'_{\mathscr X,v} .
\end{align*}
C'est un produit eulérien absolument convergent pour $\Re(s)>-\eps$,
d'où un prolongement par continuité en $s=0$, de valeur
\begin{multline*}
\prod_{v\nmid\infty} \zeta_v(1)^{d-\# J}
     \mu'_{\mathscr X,v}(\mathscr X(F_v)) 
      \prod_{v|\infty} \mu'_{\mathscr X,v} (\mathscr X(F_v)) \\
= \tau(\mathscr X) \mu(\A_F/F)^d
(\res_{s=1}\zeta_F(s))^{-\rg(\Pic\mathscr X_F)} 
\end{multline*}
en vertu de la définition~\eqref{defi:mestam} de la mesure
de Tamagawa de $\mathscr X(\A_F)$.
Ainsi,
\begin{align*}
\lim_{s\ra 0} \check H(-s\lambda-\rho,\mathbf 1) (\prod_j s\lambda_j) &
= (\res_{s=1}\zeta_F(s))^{\# J}
   \lim_{s\ra 0}
     \check H(-s\lambda-\rho,\mathbf 1)
	\prod_j \zeta_F(1+s\lambda_j)^{-1}
 \\
&= \tau(\mathscr X) \mu(\A_F/F)^d (\res_{s=1}\zeta_F(s))^d. 
\end{align*}
Finalement, on a donc
\begin{align*}
 \lim_{s\ra 0} F(\lambda s)\mathsf X_{\PL^+(\Sigma)}(\lambda s)^{-1}
   &   = (\vol(\A_F/F)\res_{s=1}\zeta_F(s))^{-d}
      \mu(\A_F/F)^d (\res_{s=1}\zeta_F(s))^d
            \tau(\mathscr X) \\
 &  = \tau(\mathscr X), 
\end{align*}
ainsi qu'il fallait démontrer.
\end{proof}

L'équation~\eqref{eq:poisson.b} et le théorème~\ref{theo:prolongement}
impliquent alors le théorème suivant.
\begin{theo}\label{theo:torique}
La fonction zêta des hauteurs (décalée)
\[ \lambda \mapsto Z(\rho+\lambda) \]
converge localement uniformément sur le tube
$\Tube(\PL(\Sigma)^+)$ et définit une fonction holomorphe
sur $\Tube(\Lambda^0_\eff(\mathscr X_F))$.
Si $\beta>1$ et si $\mathscr D_\beta$ désigne la classe de contrôle introduite
au sous-paragraphe~\ref{subsec:controle},
elle appartient à l'espace
$\mathscr H_{\{0\}}(\Lambda^0_\eff(\mathscr X_F);\Lambda^0_\eff(\mathscr
X_F))$
(défini en~\ref{defi:H(C)})
des fonctions méromorphes $\{0\}$-contrôlées dont les pôles
sont simples et donnés par les faces du cône $\Lambda^0_\eff(\mathscr X_F)$.

De plus, pour tout $\lambda\in\Lambda^0_\eff(\mathscr X_F)$,
\[ \lim_{s\ra 0} \frac{Z(s\lambda+\rho)}{\mathsf X_{\Lambda^0_\eff}(s\lambda)}
         = \tau(\mathscr X). \]
\end{theo}

En spécialisant la fonction zêta des hauteurs
à la droite $\C\rho$ qui correspond au fibré en droite anticanonique,
on obtient le corollaire:
\begin{coro}
Si $\beta>1$, il existe $\eps>0$, une fonction $f$ holomorphe
pour $\Re(s)\geq 1-\eps$ telle que
\begin{enumerate}\def\theenumi{\roman{enumi}}\def\labelenumi{(\theenumi)}
\item
$f(1)= \tau(\mathscr X)$ ;
\item
Pour tout $\sigma\in [1-\eps;1+\eps]$ et tout $\tau\in\R$,
$\abs{f(\sigma+i\tau)}\ll (1+\abs {\tau})^{\beta} $ ;
\item
Pour tout $\sigma>1$ et tout $\tau\in\R$,
$Z(s\omega) = \big(\frac{s}{s-1}\big)^r f(s)$.
\end{enumerate}
\end{coro}

\begin{coro}
Si $r$ désigne le rang de $\Pic(\mathscr X_F)$,
il existe un polynôme unitaire $P$  de degré
$r-1$ et
un réel $\eps>0$ tels que pour tout $H > 0$,
\[ N(\omega_{\mathscr X}^{-1};H)
       = \frac{\tau(\mathscr X )}{(r-1)!}
          H P(\log H)  + O(H^{1-\eps}). \]
\end{coro}
Lorsque $F=\Q$ et lorsque la variété torique $\mathscr X$
est projective et telle que $\omega_{\mathscr X}^{-1}$ est engendré
par ses sections globales, ce corollaire avait été démontré
précédemment par R.~de la Bretèche.
Sa méthode est différente ; elle est
fondée sur le travail de P.~Salberger~\cite{salberger98}
et une étude fine des sommes de fonctions arithmétiques
en plusieurs variables
(voir \cite{breteche98c,breteche98b} et \cite{breteche98}
pour un cas particulier).


\section{Application aux fibrations en vari\'et\'es toriques}
\label{sec:fibrations}%

\subsection{Holomorphie}
\label{subsec:holomorphie}%

Soit $\mathcal B$ un $S$-sch\'ema projectif et plat.
Soit $\mathcal T\ra\mathcal B$ un $G$-torseur, et notons
$\eta:X^*(G)\ra\Pic(\mathcal B)$ l'homomorphisme 
de fonctorialit\'e des torseurs.
Fixons un rel\`evement $\hat\eta:X^*(G)\ra\hatPic(\mathcal B)$
de cet homomorphisme (c'est-\`a-dire, un choix de m\'etriques hermitiennes
\`a l'infini sur les images d'une base de $X^*(G)$, prolong\'es
par multiplicativit\'e \`a l'image de $\eta$).

Donnons nous une $S$-vari\'et\'e torique lisse $\mathcal X$, compactification
\'equivariante de $G$. Soit $\mathcal Y$ le $S$-sch\'ema obtenu par
les constructions du \S\,\ref{constructions}.

On obtient alors un diagramme canonique, qui provient des
propositions~\ref{prop:iota-eta}, \ref{prop:hat-iota-eta},
du th\'eor\`eme~\ref{theo:picard} et de l'oubli des m\'etriques hermitiennes :
\begin{equation}
\xymatrix{
 0 \ar[r] & {X^*(G)} \ar[r]
    & {\Pic^G(\mathcal X_F)\oplus \Pic(\mathcal B_F)} \ar[r]
    & { \Pic(\mathcal Y_F) }  \ar[r] & 0 \\
 0 \ar[r] &  { X^*(G) } \ar[r] \ar@{=}[u] 
    & { \hatPic^{G,\mathbf K}(\mathcal X) \oplus \hatPic(\mathcal B) }
        \ar[u]  \ar[r]
    & { \hatPic(\mathcal Y) } \ar[u] 
}
\end{equation}

Le sch\'ema $\mathcal Y$ contient $\mathcal T$ comme ouvert dense.
On s'int\'eresse \`a la fonction z\^eta des hauteurs de $\mathcal T$.
Lorsque $\lambda\in\Pic^{G}(\mathcal X_F)_\C$,
notons $\hat\lambda$ l'image de $\lambda$ par
l'homomorphisme~\eqref{eq:metcan}.
Si de plus $\hat\alpha\in\hatPic(\mathcal B)$, 
on notera enfin
\[ Z(\hat\lambda,\hat\alpha) =
Z(\vartheta(\hat\lambda)\otimes\pi^*\hat\alpha; \mathcal Y)
= \sum_{y\in \mathcal T(F)}
        H( \vartheta(\hat\lambda)\otimes\pi^*\hat\alpha; y)^{-1}.  \]

\begin{prop}\label{prop:poisson}%
Soient $\widehat\Lambda\subset\hatPic (\mathcal B)_\R$ une partie
convexe telle que $Z(\hat\alpha;\mathcal B)$ converge normalement
si la partie r\'eelle de $\hat\alpha\in\hatPic(\mathcal B)_\C$
appartient \`a $\Lambda$.

Alors, la fonction z\^eta des hauteurs de $\mathcal T$ converge
absolument pour tout $(\hat\lambda,\hat\alpha)$ tel que la partie r\'eelle de 
$\lambda\otimes\omega_{\mathcal X}$ appartient \`a
$\Lambda_{\eff}^0 (\mathcal X_F)$
et la partie r\'eelle de $\alpha$ appartient \`a $\Lambda$.
La convergence est de plus uniforme si la partie la partie
r\'eelle de $\lambda\otimes\omega_{\mathcal X}$ d\'ecrit un compact
de $\Lambda_{\eff}^0 (\mathcal X_F)$.
\end{prop}
\begin{proof}
On peut d\'ecomposer la fonction z\^eta des hauteurs de $\mathcal T$ en
\'ecrivant
\begin{equation}\label{eq:decomposition}%
Z(\hat\lambda,\hat\alpha)  = \sum_{b\in \mathcal B(F)}
      H(\hat\alpha;b)^{-1} Z (\vartheta(\hat\lambda); \mathcal T|_b) . 
\end{equation}
D'apr\`es la remarque~\ref{rema:torsion-torique},
le fibr\'e inversible $\lambda$ admet une section $G$-invariante $\mathsf s$
qui n'a ni p\^oles ni z\'eros sur l'ouvert $G\subset\mathcal X$.
En utilisant cette section, on obtient, en vertu
du th\'eor\`eme~\ref{theo:comp-torsarith} et
de la proposition~\ref{prop:comp-tordue}
une \'egalit\'e
\begin{equation}
Z(\vartheta(\hat\lambda);\mathcal T|_b)
      = \sum_{x\in G(F)} H(\hat\lambda,\mathsf s;\mathbf g_b \cdot x)^{-1},
\end{equation} 
o\`u $\mathbf g_b\in G(\A_F)$ repr\'esente la classe du $G$-torseur
arithm\'etique $\widehat{\mathcal T}|_b$. 
On rappelle que si $\mathbf x\in G(\A_F)$, on a une expression
de la hauteur en produit de hauteurs locales
\[ H(\hat\lambda,\mathsf s,\mathbf x)
    = \prod_{v} \norm{\mathsf s}_v(x_v)^{-1} . \]
On peut appliquer la formule sommatoire de Poisson sur le tore ad\'elique
$G(\A_F)$, d'o\`u, en utilisant l'invariance des hauteurs locales
par les sous-groupes compacts maximaux,
\begin{equation}\label{eq:poisson}%
 Z(\vartheta(\hat\lambda);\mathcal T|_b) =
  \int_{\mathcal A_G}
      \chi^{-1}(\mathbf g_b) \Fourier H (-\hat\lambda;\chi)\, d\chi 
\end{equation}
o\`u l'int\'egration est sur le groupe $\mathcal A_G$
des caract\`eres (unitaires continus)
du groupe localement compact $G(F)\backslash G(\A_F)/K_G$, muni de son unique mesure de Haar $d\chi$
qui permet cette formule.

L'utilisation de la formule de Poisson est justifi\'ee par le fait
que les deux membres convergent absolument.
La s\'erie du membre de gauche est trait\'ee
dans~\cite{batyrev-t98b}, Theorem 4.2,
lorsque $\mathbf g_b=1$, c'est-\`a-dire lorsqu'il n'y a pas de torsion.
Comme il existe une constante $C(\lambda,\mathbf g_b)$ ne d\'ependant que
de $\mathbf g_b$ et $\hat\lambda$ telle que
\[
\abs{ H(\hat\lambda,\mathsf s; \mathbf g_b\cdot x ) }^{-1}
\leq C(\lambda,\mathbf g_b) \abs{ H(\hat\lambda,\mathsf s; x) } ^{-1}
\]
et comme $ H(\hat\lambda,\mathsf s; x ) = H(\hat\lambda;x) $,
la convergence absolue
du membre de gauche en r\'esulte.
(Voir aussi le lemme~\ref{lemm:convergence}.)
Quant \`a l'int\'egrale du membre de droite, on peut n\'egliger le
caract\`ere $\chi$ dont la valeur absolue est~$1$ et on retrouve
une int\'egrale dont la convergence absolue est prouv\'ee
dans~\cite{batyrev-t98b} (preuve du th\'eor\`eme 4.4).
Cela prouve aussi que lorsque $\Re(\lambda)$ d\'ecrit un compact
de $\omega_{\mathcal X}^{-1} + \Lambda_{\eff}^0(\mathcal X_F)$,
la fonction z\^eta des hauteurs $Z(\vartheta(\hat\lambda);\mathcal T|_b)$
de la fibre en $b\in\mathcal B(F)$
est born\'ee ind\'ependamment de $b$.

En reportant cette majoration dans la
d\'ecomposition~\eqref{eq:decomposition}, il en r\'esulte la convergence
absolue de la fonction z\^eta des hauteurs de $\mathcal T$ lorsque
la partie r\'eelle de $\hat\alpha$ appartient \`a $\widehat\Lambda$
et $\lambda\otimes\omega_{\mathcal X}$ appartient \`a
$\Lambda_{\eff}^0(\mathcal X_F)$, uniform\'ement
lorsque $\lambda\otimes\omega_{\mathcal X}$ d\'ecrit un compact de ce c\^one.
\end{proof}

Dans~\cite{chambert-loir-t99}, définition~\ref{defi:LArakelov},
on a défini la notion de fonction~$L$ d'Arakelov
attachée à un torseur arithmétique et à une fonction sur un espace
adélique.
Appliquée au $G\times\gm$-torseur arithm\'etique
sur $\mathcal B$ d\'efini par
$\widehat{\mathcal T}\times_{\mathcal B} \hat\alpha$
et à la fonction $\chi^{-1} \cdot \norm{\cdot}$,
la définition devient
\[ L ( \widehat {\mathcal T}\sqtimes \hat\alpha, \chi^{-1}\sqtimes
\norm{\cdot})
=  \sum_{b\in\mathcal B(F)} \chi^{-1}(\mathbf g_b) H(\hat\alpha;b)^{-1} . \]
(On a utilisé le fait que $\mathbf g_b\in G(F)\backslash G(\A_F)/K_G$
est la classe du $G$-torseur arithm\'etique $\mathcal T|_b$.)

Un corollaire de la d\'emonstration de la proposition précédente
est alors le suivant :
\begin{coro}
Sous les hypoth\`eses de la proposition~\ref{prop:poisson}, on a la formule
\[ 
Z(\hat\lambda,\hat\alpha) = \int_{\mathcal A_G} \Fourier H(-\hat\lambda;\chi)
L ( \widehat {\mathcal T}\sqtimes  \hat\alpha, \chi^{-1}\sqtimes \norm{\cdot})
\, d\chi.
\] 
\end{coro}
\begin{proof}
Compte tenu de la majoration \'etablie \`a la fin de la preuve du th\'eor\`eme
pr\'ec\'edent et des rappels faits sur les fonctions $L$ d'Arakelov,
il suffit de reporter l'\'equation~\eqref{eq:poisson} dans la
formule~\eqref{eq:decomposition} et d'\'echanger les signes somme et
int\'egrale.
\end{proof}

Cette derni\`ere formule est le point de d\'epart pour \'etablir,
moyennant des
hypoth\`eses suppl\'ementaires sur $\mathcal B$,
un prolongement m\'eromorphe de la fonction z\^eta des hauteurs
de $\mathcal T$.

\subsection{Prolongement m\'eromorphe}
\label{subsec:montee}%

Fixons une section de l'homomorphisme canonique
$\hatPic(\mathcal B)\otimes_\Z\Q \ra \Pic(\mathcal B_ F)\otimes\Q$,
autrement dit un choix de fonctions hauteurs compatible
au produit tensoriel, ce que Peyre appelle \emph{syst\`eme de hauteurs}
dans~\cite{peyre98}, 2.2.
Concernant $\mathcal X$, on utilise toujours les m\'etriques
ad\'eliques canoniques utilis\'ees au paragraphe~\ref{sec:torique}.
Ainsi, on \'ecrira $\lambda$ et $\alpha$, les chapeaux devenant inutiles.
L'application $\hat\eta:X^*(G)\ra \hatPic(\mathcal B)$
est suppos\'ee \^etre la compos\'ee de l'application
$\eta:X^*(G)\ra\Pic(\mathcal B_F)$ donn\'ee
par la restriction du torseur \`a la fibre g\'en\'erique,
et de la section
$\Pic(\mathcal B_F)\otimes\Q\ra\hatPic(\mathcal B)\otimes\Q$ fix\'ee.

Ces restrictions ne sont pas vraiment essentielles mais simplifient
beaucoup les notations.

Notons $V_1=\Pic^G(\mathcal X_ F)_\R$,
$M_1=X^*(G)_\R$,
$n_1=\dim V_1$ et
$V_2=\Pic(\mathcal B_ F)_\R$.
Soient $\Lambda_1\subset V_1$ et $\Lambda_2\subset V_2$ les c\^ones ouverts,
int\'erieurs des c\^ones effectifs dans $\Pic^G(\mathcal X_F)_\R$
et $\Pic(\mathcal B_ F)_\R$.
L'espace vectoriel $V_1$ poss\`ede une base naturelle, form\'ee
des fibr\'es en droites $G$-lin\'earis\'es associ\'es aux
diviseurs $G$-invariants sur $\mathcal X\otimes F$. Dans
cette base, le c\^one $\Lambda_1$ est simplement l'ensemble
des $(s_1,\ldots,s_{n_1})$ strictement positifs.

On note $\eta:M_1\ra V_2$ l'application lin\'eaire 
d\'eduite de $\hat\eta$ et $M=(\id,-\eta)(M_1)\subset V_1\times V_2$.
Notons $V=V_1\times V_2$. Les th\'eor\`emes~\ref{theo:picard}
et~\ref{theo:effectif} identifient $\Pic(\mathcal Y_F)_\R$ \`a
$V/M$, et l'int\'erieur du c\^one effectif de $\mathcal Y_F$
\`a l'image de $\Lambda_1\times \Lambda_2$ par la projection $V\ra V/M$.
Si $\omega_{\mathcal X}$ est muni de sa $G$-lin\'earisation canonique,
la proposition~\ref{prop:canonique} dit que $\omega_{\mathcal Y}$
est l'image du couple $(\omega_{\mathcal X},\omega_{\mathcal B})$
par cette m\^eme projection.

\begin{lemm}
Avec ces notations,
la formule du corollaire du paragraphe pr\'ec\'edent peut 
se r\'ecrire :
\[ Z(\lambda+\omega_{\mathcal X}^{-1},\alpha+\omega_{\mathcal B}^{-1})
 = \int_{M_1} f(\lambda+im_1;\alpha-i\eta(m_1))\, dm_1, \]
o\`u la fonction
\[ f:\Tube(\Lambda_1 \times\Lambda_2)\ra\C \]
est d\'efinie par
\[ f(\lambda;\alpha)=
\int_{\mathcal U_G} 
  \Fourier H(-(\lambda+\omega_{\mathcal X}^{-1});\chi_u)
  L (\widehat{\mathcal T}\sqtimes(\alpha+\omega_{\mathcal B}^{-1});
        \chi_u^{-1}\sqtimes \norm{\cdot})
  \, d\chi_u
\]
et $dm_1$, $d\chi_u$ sont des mesures de Haar sur $M_1$
et $\mathcal U_G$ telles que $d\chi=dm_1\,d\chi_u$ dans
la d\'ecomposition $\mathcal A_G=M_1\oplus\mathcal U_G$
du paragraphe~\ref{rappel:caracteres} (cf.~aussi le
lemme~\ref{lemm:normalisations}).
\end{lemm}
On note que $\mathcal U_G$ est un groupe discret et que
la mesure $d\chi_u$ est donc proportionnelle \`a la mesure
de comptage.
\begin{proof}
Si $\chi\in\mathcal A_G$ s'\'ecrit
$(m_1,\chi_u)$ dans $ M_1\oplus \mathcal U_G$,
on remarque que l'on a les \'egalit\'es
\[ \Fourier H(-\lambda;\chi)= \Fourier H(-\lambda-\iota(i\,m_1);\chi_u)\]
et 
\[ \chi^{-1}(\mathbf g_b) H(\hat \alpha ;b)^{-1}
 = \chi_u^{-1}(\mathbf g_b) H(\hat \alpha-\eta(m_1);b)^{-1} \]
car (lemme~\ref{lemm:chi/H})
\[ \chi_{m_1}(\mathbf g_b) = \exp (i\norm{\cdot}) ([m_1]_* \widehat
{\mathcal T}|_b)
=\exp(i\norm{\hat\eta(m_1)|_b})
=H(-\hat\eta(m_1);b).
\]
On utilise ensuite le th\'eor\`eme de Fubini.
\end{proof}

On utilise enfin les notations du \S\,\ref{sec:analyse}.

\begin{enonce}{Hypoth\`eses}
On fait les hypoth\`eses suivantes : 
\begin{itemize}
\item le cône $\Lambda_2$ est un cône poly\'edral (de type fini). Notons
$(\ell_j)$ les formes linéaires définissant ses faces ;
\item la fonction z\^eta des hauteurs de $\mathcal B$ converge
localement normalement pour $\alpha+\omega_{\mathcal B}\in \Lambda_2$ ;
\item
il existe un voisinage convexe $B_2$ de l'origine dans $V_2$
et
pour tout caract\`ere $\chi\in\mathcal A_G$
une fonction holomorphe $g(\chi;\cdot)$ sur le tube
$\Tube(B_2)$ tels que, si $\Re(\alpha+\omega_{\mathcal B})\in\Lambda_2$,
\[ L ( \widehat{\mathcal T}\sqtimes \alpha,\chi^{-1} \sqtimes \norm{\cdot})
 = \prod_j \frac{ \ell_j(\alpha)}{\ell_j(\alpha+\omega_{\mathscr B})}
   g(\chi;\alpha+\omega_{\mathcal B}) ; \]
\item
il existe un réel $\gamma$ strictement positif tel que pour tout $\eps>0$,
les fonctions $g(\chi;\cdot)$ v\'erifient 
une majoration uniforme
\[ \abs{g(\chi;\alpha+\omega_{\mathcal B})}
          \leq C_\eps \big(1+\norm{\Im(\alpha)}\big)^{\gamma}
                      \big(1+ \norm{\chi}\big)^{\eps}, \]
pour un r\'eel $\eps<1$ et une constante $C_\eps$ ;
\item
si $\tau(\mathscr B)$ désigne le nombre de Tamagawa de $\mathscr B$,
pour tout $\alpha$ appartenant à $\Lambda_2$,
\[ \lim_{s\ra 0^+} \frac{Z(\mathscr B; s\alpha+\omega_{\mathscr B}^{-1})}
           {\mathsf X_{\Lambda_2}(s\alpha)}
     = \tau(\mathscr B) \neq 0 .\]
\end{itemize}
\end{enonce}
\begin{rema}
Dans le cas où $\mathscr B$ est une variété de drapeaux généralisée,
ces hypothèses correspondent à des énoncés sur les séries d'Eisenstein
tordues par des caractères de Hecke. Ils sont établis
dans~\cite{strauch-t99}.
\end{rema}

Dans la suite, on travaille
avec les classes de contr\^ole $\mathscr D_\beta$
introduites au  paragraphe~\ref{subsec:controle}.

\begin{lemm}
Sous les hypothèses précédentes, pour tout réel $\beta>1$,
la fonction $f$ appartient \`a $\mathcal H_M(\Lambda_1\times \Lambda_2)$,
pour la classe $\mathscr D_{\beta+\gamma}$.
\end{lemm}
\begin{proof}
Il suffit de reprendre la démonstration
de la proposition~\ref{prop:Fcontrole}, d'y insérer les majorations
que nous avons supposées
et de majorer
\[ (1+\norm{\Im\lambda})^\beta (1+\norm{\Im\alpha})^\gamma 
 \leq (1+\norm{\Im\lambda}+\norm{\Im\alpha})^{\beta+\gamma}. \]
\end{proof}

Grâce au théorème d'analyse~\ref{theo:prolongement},
on en déduit un prolongement méromorphe pour la fonction zêta
des hauteurs de $\mathscr T$.
\begin{theo}\label{theo:montee}
La fonction zêta des hauteurs décalée de $\mathscr T$ 
admet un prolongement méromorphe dans un voisinage de
$\Tube(\Lambda^0_\eff(\mathscr Y))$
dans $\Pic(\mathscr Y)_\C$.
Cette fonction a des pôles simples
donnés par les équations des faces de $\Lambda^0_\eff(\mathscr Y)$.
De plus, si $\lambda\in\Lambda^0_\eff(\mathscr Y)$,
\[ \lim_{s\ra 0^+} \frac{Z(\mathscr T;s\lambda+\omega_{\mathscr Y}^{-1})}
{\mathsf X_{\Lambda_\eff(\mathscr Y)}(s\lambda)}
         = \tau(\mathscr Y), \]
le nombre de Tamagawa de $\mathscr Y$.
\end{theo}
\begin{proof}
Le seul point qui n'a pas été rappelé est que le nombre
de Tamagawa est $\mathscr Y$ est le produit de ceux de
$\mathscr X$ et $\mathscr B$ (\cite{chambert-loir-t99}, théorème 2.5.5).
\end{proof}

\begin{coro}
Il existe $\eps>0$ et un polynôme $P$ tels que
le nombre de points de $\mathscr T(F)$ dont la hauteur  anticanonique est
inférieure ou égale à $H$ vérifie un développement asymptotique
\[ N(H) = H P(\log H) + O(H^{1-\eps}) \]
lorsque $H$ tend vers~$+\infty$.
Le degré de $P$ est égal au rang de $\Pic(\mathscr Y_F)$
moins~$1$
et son coefficient dominant vaut
\[ \mathsf X_{\Lambda_\eff(\mathscr Y)}(\omega_{\mathscr Y}^{-1})
         \tau(\mathscr Y). \]
\end{coro}

\bgroup \small
\appendix
\def\thesection{\Alph{section}}
\def\theparagraph{\thesection.\arabic{paragraph}}
\let\theequation\theparagraph

\vskip 0pt plus .3\textheight
\penalty-100
\vskip 0pt plus -.3\textheight

\setcounter{paragraph}{0}

\section{Un th\'eor\`eme taub\'erien}

Le but de ce paragraphe est de d\'emontrer un th\'eor\`eme taub\'erien
dont la preuve nous a été communiquée par P.~Etingof.
Ce théorème est certainement bien connu des experts
mais que nous n'avons pu le trouver sous cette forme dans la littérature.

\begin{theo}\label{theo:tauber}%
Soient $(\lambda_n)_{n\in\N}$ une suite croissante
de r\'eels strictement positifs, $(c_n)_{n\in\N}$ une suite
de réels positifs et $f$ la s\'erie de Dirichlet
\[ f(s) = \sum_{n=0}^\infty  c_n \frac{1}{\lambda_n^s}. \]
On fait les hypoth\`eses suivantes :
\begin{itemize}
\item la s\'erie d\'efinissant $f$ converge dans un demi-plan $\Re(s)>a>0$ ;
\item elle admet un prolongement m\'eromorphe dans un demi-plan
$\Re(s)> a-\delta_0>0$ ;
\item dans ce domaine, elle poss\`ede un unique p\^ole en $s=a$,
de multiplicit\'e $b\in\N$. On note $\Theta=\lim_{s\ra a} f(s)(s-a)^b>0$ ;
\item  enfin, il existe un r\'eel $\kappa>0$ de sorte
que l'on ait pour $\Re(s)>a-\delta_0$ l'estimation, 
\[ \abs{f(s) \frac{(s-a)^b}{s^b}} = O \big((1+\Im(s))^{\kappa}\big). \]
\end{itemize}
Alors il existe un polyn\^ome unitaire $P$ de degr\'e $b-1$ tel que pour tout
$\delta<\delta_0$, on ait, lorsque $X$ tend vers $+\infty$,
\[ N(X) \stackrel{\text{def}}= \sum_{\lambda_n\leq X} c_n
    = \frac{\Theta}{a\, (b-1)!} X^a P(\log X) + O(X^{a-\delta}). \]
\end{theo}

On introduit pour tout entier $k\geq 0$ la fonction
\[ \phi_k (X) = \sum_{\lambda_n\leq X} a_n \left(\log (X/\lambda_n)\right)^k, \]
de sorte que $\phi_0=N$.
\begin{lemm}
Sous les hypoth\`eses du th\'eor\`eme~\ref{theo:tauber},
il existe pour tout entier $k>\kappa$
un polyn\^ome $Q$ de degr\'e $b-1$ et de coefficient
dominant $k!\Theta /(a^{k+1}(b-1)!)$
tel que pour tout $\delta<\delta_0$,
on ait l'estimation, lorsque $X$ tend vers $+\infty$,
\[ \phi_k(X)
    =  X^a Q(\log X) + O(X^{a-\delta}). \]
\end{lemm}

\begin{proof}
Soit $a'>a$ arbitraire.
On remarque, en vertu de l'int\'egrale classique
\[  \int_{a'+i\R}  \lambda^s \frac{ds}{s^{k+1}}
    = \frac{2i\pi}{k!}\left(\log^+(\lambda)\right)^k, \quad \lambda >0 \]
que l'on a la formule
\[ \phi(X) = \frac{k!}{2i \pi} \int_{a'+i\R} f(s) X^s \frac{ds}{s^{k+1}}, \]
l'int\'egrale \'etant absolument convergente puisque $\kappa<k$.

On veut d\'ecaler le coutour d'int\'egration vers la droite
verticale $\Re(s)=a-\delta$, o\`u $\delta$ est un r\'eel arbitraire
tel que $0<\delta<\delta_0$.
Dans le rectangle $a-\delta\leq\Re(s)\leq a'$, $\abs{\Im(s)}\leq T$,
il y a un unique p\^ole en $s=a$. Le r\'esidu y vaut
\[ \Res_{s=a} f(s) \frac{X^s}{s^{k+1}} =
    \frac{\Theta}{a^{k+1}\, (b-1)!} X^a Q(\log X) \]
o\`u $Q$ est un polyn\^ome unitaire de degr\'e $b-1$.
Il en r\'esulte que
\begin{multline*}
 \frac{1}{2i\pi} \int_{a'-iT} ^{a'+iT} f(s) X^s\frac{ds}{s^{k+1}} \\
     = \frac{1}{2i\pi} \int_{a-\delta-iT}^{a-\delta+iT} f(s)
X^s\frac{ds}{s^{k+1}} + I_+ - I_-
        + \frac{\Theta}{a^{k+1}\, (b-1)!} X^a Q(\log X) ,
\end{multline*}
o\`u $I_+$ et $I_-$ sont les int\'egrales sur les segments horizontaux
(orient\'es de la gauche vers la droite).
Lorsque $T$ tend vers $+\infty$, ces int\'egrales sont 
$O(T^{\kappa-k-1}X^{a'})$
et tendent donc vers~$0$.
Les hypoth\`eses sur $f$ et le fait que $k>\kappa$
montrent que $f(s)X^s/s^{k+1}$ est absolument
int\'egrable sur la droite $\Re(s)=a-\delta$, l'int\'egrale \'etant
major\'ee par $O(X^{a-\delta})$.
Par cons\'equent, on a
\[ \phi(X) = \Theta \frac{k!}{a^{k+1}\, (b-1)!} X^a Q(\log X) 
    + O(X^{a-\delta}).
\]
Le lemme est ainsi d\'emontr\'e.
\end{proof}

\begin{proof}[Preuve du th\'eor\`eme]
On va démontrer par récurrence descendante que la conclusion du lemme
précédent vaut en fait pour tout entier $k\geq 0$. Arrivés à $k=0$, le
théorème sera prouvé. Montrons donc comment passer de
$k\geq 1$ à $k-1$.

Pour tout $\eta\in\left]0;1\right[$, on a facilement l'in\'egalit\'e
\[
\frac{\phi_k(X(1-\eta))-\phi_k(X)}{\log(1-\eta)}
\leq k \phi_{k-1}(X) \leq
\frac{\phi_k(X(1+\eta))-\phi_k(X)}{\log(1+\eta)} . \]
Fixons un r\'eel $\delta'$ tel que $0<\delta'<\delta<\delta_0$.
D'apr\`es le lemme pr\'ec\'edent, il existe un r\'eel $C$ tel que
\[ \abs{\phi_k(X) - \frac{k!\Theta}{a^{k+1}\, (b-1)!} X^a Q(\log X) } \leq C
X^{a-\delta'}. \]
On constate que l'on a alors, si $-1<u<1$,
\begin{multline*} \frac{\phi_k(X(1+u))-\phi_k(X)}{\log(1+u)}  \\
  =  \frac{k!\Theta }{a^{k+1}\, (b-1)!} X^a \frac{Q(\log X+\log(1+u))(1+u)^a
         - Q(\log X)}{\log(1+u)} + R(X), 
\end{multline*}
o\`u
\[ \abs{R(X)} \leq 2C X^{a-\delta'} /\abs{\log(1+u)}
           = O(X^{a-\delta'}/u)  \]
si $u$ tend vers~$0$ et $X\ra+\infty$.
Toujours lorsque $X\ra +\infty$ et $u\ra 0$, on a
\begin{align*}
\frac{Q(\log X+\log(1+u))(1+u)^a - Q(\log X)}{\log(1+u)} \hskip -15em \\
& = Q(\log X) \frac{(1+u)^a-1}{\log(1+u)} + \sum_{n=1}^{b-1} 
\frac{1}{n!} Q^{(n)}(\log X) \log(1+u)^{n-1} (1+u)^a \\
&= Q (\log X) \left (a + O(u) \right) 
 + Q'(\log X) \left (1 + O(u) \right)
 +   O((\log X)^{b-1} u) \\
&= (a Q+Q')(\log X)+  O ((\log X)^{b-1} u) .
\end{align*}
Prenons $u=\pm 1/X^\eps$ où $\eps>0$ est choisi de
sorte que $\delta'+\eps<\delta$.
Alors, lorsque $X\ra+\infty$, $\abs{R(X)} =O( X^{a-\delta}) $
et
\[ \frac{Q(\log X+\log(1+u))(1+u)^a - Q(\log X)}{\log(1+u)}
= (aQ+Q')(\log X) + O(X^{-\delta}). \]
On a alors un développement
\[
 \phi_{k-1}(X)= \frac 1k   X^a (aQ+Q')(\log X) + O(X^{a-\delta})
\]
Le coefficient dominant de $(aQ+Q')/k$ est \'egal \`a
$(k-1)!\Theta/(a^k(b-1)!)$
d'o\`u le th\'eor\`eme par récurrence descendante.
\end{proof}

\setcounter{paragraph}{0}

\section{Démonstration de quelques inégalités}
\label{app:controle}%

Le but de cet appendice est de démontrer les inégalités
sous-jacentes à la proposition~\ref{prop:controle}
qui affirmait l'existence d'une classe de contrôle.

Rappelons les notations.

Soit $\beta$ un réel strictement positif.
Si $M$ et $V$ sont deux $\R$-espaces vectoriels de dimension finie
avec $M\subset V$, notons $\mathscr D_ {\beta,\eps}(M,V)$ le sous-monoïde
de $\mathscr F(V,\R_+)$ engendré 
par les fonctions $h:V\ra\R_+$ telles qu'il existe $c>0$
et une famille $(\ell_j)$ de formes linéaires
sur $V$ vérifiant :
\begin{itemize}
\item la famille $(\ell_j|_M)$ forme une base de $M^*$ ;
\item pour tout $v\in V$ et tout $m\in M$, on a
\begin{equation}
h(v+m) \leq c \frac{(1+\norm v)^{\beta}}{(1+\norm{m})^{1-\eps}}
        \frac{1}{\prod(1+\abs{\ell_j(v+m)})}.
\end{equation}
\end{itemize}
On définit ensuite $\mathscr D_\beta(M,V)=\bigcap_{\eps>0}\mathscr
D_{\beta,\eps}(M,V)$.
\begin{theo}
Les $\mathscr D_\beta({M,V})$ définissent une classe de contrôle au
sens de la définition~\ref{defi:classe}.
\end{theo}

\begin{proof}
Les points~(\ref{defi:classe},\ref{axiome:hypomonoide})
et~(\ref{defi:classe},\ref{axiome:limite})
sont clairs. L'axiome~(\ref{defi:classe},\ref{axiome:projecteur}) est vrai
car la famille $(\ell_j\circ p|_M)$ contient une base de $(M/M_1)^*$.
L'axiome~(\ref{defi:classe},\ref{axiome:max}) résulte de l'inégalité
\[ \min_{\abs{t}\leq 1} (1+\abs{\ell(v+tu+m)}) \geq
\frac{1}{1+\abs{\ell(u)}} (1+\abs{\ell(v+m)} \]
valable pour tous $v\in V$, $u\in V$ et $m\in M$.
Enfin, l'axiome~(\ref{defi:classe},\ref{axiome:integrale}), le plus
délicat, fait l'objet de la proposition suivante.
\end{proof}

\begin{prop}\label{prop:majoration}
Soient $M\subset V$, $V'$ un supplémentaire de $M$ dans $V$,
$dm$ une mesure de Lebesgue sur $M$,
$(\ell_j)$ une base de $V^*$.
Pour tout $\eps'>\eps$, il existe une constante $c_{\eps'}$ et
un ensemble $((\ell_{j,\alpha})_j)_\alpha$ de bases de $(V')^*$
tels que pour tous $v_1$ et $v_2\in M'$,
\[
\int_M \frac{1}{(1+\norm{v_1+m})^{1-\eps}}
\frac{dm}{\prod(1+\abs{\ell_j(v_2+m)})}
\leq \frac{ c_{\eps'}}{(1+\norm{v_1})^{1-\eps'}}
      \sum_\alpha \frac{1}{\prod_j (1+\abs{\ell_{j,\alpha}(v_2)})}. \]
\end{prop}

\begin{proof}
On raisonne par récurrence sur $\dim M$.
Soient $\mathbf u\in M$, $M'\subset M$ tels que $M=M'\oplus \R\mathbf u$
et fixons une mesure de Lebesgue $dm'$ sur $M'$ telle que $dm'\cdot dt=dm$.
Alors,

\begin{align*}
\int_{\R\mathbf u} \cdots
   & \ll \int_\R  \frac{1}{(1+\norm{v_1+m'}+\abs t)^{1-\eps}}
       \frac{dt}{\prod_j (1+\abs{\ell_j(v_2+m') t\ell_j(\mathbf u)})} \, dt  \\
& \ll \prod_{j\,;\,\ell_j(\mathbf u)=0}
        \frac{1}{1+\abs{\ell_j(v_2+m')}}  \times \\
& \qquad\qquad
    \times \int_\R \frac{1}{(1+\norm{v_1+m'}+\abs t)^{1-\eps}}
        \prod_{j\,;\, \ell_j(\mathbf u)\neq 0}
         \frac{1}{1+\abs{\ell_j(v_2+m')+t}}\, dt  \\
\noalign{et, en appliquant le lemme~\ref{lemm:omega} ci-dessous,} 
& \ll \frac{1+\log(1+\norm{v_1+m'})}{(1+\norm{v_1+m'})^{1-\eps}}
      \sum_\alpha \prod_j \frac{1}{1+\abs{\ell_{j,\alpha}(v_2+m')}} \\
& \ll_{\eps'} \frac{1}{(1+\norm{v_1+m'})^{1-\eps'}}
              \sum_\alpha \frac{1}{\prod_j (1+\abs{\ell_{j,\alpha}(v_2)})}.
\qquad\qed
\end{align*}
\let\qedsymbol\relax
\end{proof}
 
\begin{lemm}\label{lemm:omega}%
On a une majoration, valable pour tous réels $t_1\leq\cdots\leq t_n$
et tout $A\geq 0$,
\[ \int_{-\infty}^\infty \frac{1}{(1+A+\abs t)^{1-\eps}} \prod_{j=1}^n
\frac{1}{1+\abs{t-t_j}}\, dt 
\ll \frac{1+\log(1+A)}{(1+A)^{1-\eps}}
\sum_{\alpha}\prod_{j=1}^{n-1} \frac{1}{1+\abs{\tau_{\alpha,j}}}
\]
où pour tout $\alpha$ et tout $j$,
$\tau_{\alpha,j}=t_{a(\alpha,j)}-t_{b(\alpha,j)}$
de sorte que pour tout $\alpha$,
notant $(e_1,\ldots,e_n)$ la base canonique de $\R^n$,
les familles $(e_{\alpha,j}=e_{a(\alpha,j)}-e_{b(\alpha,j)})_j$
sont libres.
\end{lemm}
\begin{proof}
On découpe l'intégrale en $\int_{-\infty}^{t_1}$, $\int_{t_1}^{t_2}$,
\dots, $\int_{t_n}^\infty$ et on majore chaque terme.

Pour l'intégrale de $-\infty$ à $t_1$, on a
\begin{align*}
\int_{-\infty}^{t_1} \cdots & \leq \prod_{j=2}^n \frac{1}{1+\abs{t_j-t_1}}
    \int_{-\infty}^{t_1}  \frac{1}{(1+A+\abs t)^{1-\eps}} \frac{dt}{1+t_1-t}  \\
&\leq \prod_{j=2}^n \frac{1}{1+\abs{t_j-t_1}}
     \int_0^\infty \frac{1}{(1+A+\abs{t-t_1})^{1-\eps}} \frac{dt}{1+t} \\
& \leq \prod_{j=2}^n \frac{1}{1+\abs{t_j-t_1}} \frac{1+\log(1+A)}{(1+A)^{1-\eps}}
\end{align*}
d'après le lemme~\ref{lemm:alpha}.
La dernière intégrale (de $t_n$ à $+\infty$) se traite de même.
Enfin,
\begin{multline*}
\int_{t_k}^{t_{k+1}} \cdots \leq
    \prod_{j<k} \frac{1}{1+\abs{t_k-t_j}} 
    \prod_{j>k+1} \frac{1}{1+\abs{t_{k+1}-t_j}} \times \\
    \times 
    \int_{t_k}^{t_{k+1}} \frac{1}{(1+A+\abs t)^{1-\eps}}
\frac{dt}{(1+t-t_k)(1+t_{k+1}-t)} 
\end{multline*}
et cette dernière intégrale s'estime comme suit :
\begin{align*}
\hbox to 0pt{$\displaystyle
    \frac{1}{(1+A+\abs t)^{1-\eps}} \frac{dt}{(1+t-t_k)(1+t_{k+1}-t)}
={}$\hss}\hskip 2cm
 \\
&= \int_{t_k}^{t_{k+1}} \frac{1}{(1+A+\abs t)^{1-\eps}}
        \frac{1}{2+t_{k+1}-t_k} \left(
\frac{1}{1+t-t_k}+\frac1{1+t_{k+1}-t}\right) \, dt 
\\
&\leq \frac{1}{2+t_{k+1}-t_k} \Big(
      \int_{t_k}^\infty \frac{1}{(1+A+\abs t)^{1-\eps}}\frac{dt}{1+t-t_k}\\
    &\qquad\qquad
            {} + \int_{-\infty}^{t_{k+1}} \frac{1}{(1+A+\abs t)^{1-\eps}}
            \frac{dt}{1+t_{k+1}-t} \Big)
\\
&\leq \frac{1}{2+t_{k+1}-t_k} \Big(
      \int_{0}^\infty \frac{1}{(1+A+\abs{t+t_k})^{1-\eps}}\frac{dt}{1+t}
            \\
  &   \qquad\qquad
    {} + \int_{0}^{\infty} \frac{1}{(1+A+\abs {t-t_{k+1}})^{1-\eps}}
        \frac{dt}{1+t}
                            \Big) 
\\
&\ll \frac{1}{1+t_{k+1}-t_k} \frac{1+\log(1+A)}{(1+A)^{1-\eps}}
\end{align*}
en vertu du lemme~\ref{lemm:alpha}.
\end{proof}

\begin{lemm}\label{lemm:alpha}%
On a une majoration, valable pour tout $A\geq 1$ et tout $a>0$,
\[ \int_0^\infty \frac{1}{(A+\abs{t+a})^\alpha} \frac{dt}{1+t}
         \ll \frac{1+\log A}{A^\alpha}. 
\]
\end{lemm}
Il reste à démontrer ce lemme.
Pour cela, on a besoin de deux lemmes supplémentaires !
\begin{lemm}\label{lemm:+}%
Pour tous $A$ et $B\geq 1$ et tous $\alpha,\beta>0$ tels que
$\alpha+\beta>1$,
\[ \int_0^\infty \frac{dt}{(A+t)^\alpha(B+t)^\beta} \ll_{\alpha,\beta}
\frac{\min(A,B)}{A^\alpha B^\beta} \times
  \begin{cases} 1+\log(B/A) & \text{si $\alpha=1$ et $B>A$ ;} \\
                1+\log(A/B) & \text{si $\beta=1$ et $A>B$ ;} \\
                1 &\text{sinon.}
    \end{cases}
\]
\end{lemm}
\begin{proof}
On ne traite que le cas $A< B$, l'autre étant symétrique et le cas $A=B$
élémentaire.
Faisons le changement de variables $A+T=(B-A)e^u$, d'où $B+T=(B-A)(1+e^u)$.
Pour $t=0$, $u=\log A/(B-A)$. Lorsque $t\ra+\infty$, $u\ra+\infty$.
Ainsi, l'intégrale vaut
\[ I(A,B;\alpha,\beta) = \frac{1}{(B-A)^{\alpha+\beta-1}}
     \int_{\log A/(B-A)}^\infty \frac{e^{(1-\alpha)u}}{(1+e^u)^\beta} du .\]

Si $A< B\leq 2A$, on majore l'intégrale par
\begin{align*}
I(A,B;\alpha,\beta) 
&\leq \frac{1}{(B-A)^{\alpha+\beta-1}} \int_{\log A/(B-A)}^\infty
e^{(1-\alpha-\beta)u}\, du \\
&\leq \frac{1}{(B-A)^{\alpha+\beta-1}} \frac{1}{1-\alpha-\beta} \big(
\frac{B-A}A \big)^{\alpha+\beta-1} \\
&\ll \frac{1}{A^{\alpha+\beta-1}} 
\ll \frac{A}{A^\alpha B^\beta}
\end{align*}
puisque $1/A\leq 2/B$.

Lorsque $B\geq 2A$, $\log A/(B-A)\leq 0$. On minore $1+e^u$
par $1$ lorsque $u\leq 0$ et par $e^u$ lorsque $u\geq 0$, d'où
les inégalités
\begin {align*} 
(B-A)^{\alpha+\beta-1} I(A,B;\alpha,\beta)
& = 
\int_{\log A/(B-A)}^0 +  \int_0^\infty \\
& \leq  \int_0^\infty 
\frac{e^{(1-\alpha)u}}{(1+e^u)^\beta}\, du
 + \int_{\log A/(B-A)}^0 e^{(1-\alpha)u}\, du  \\
&\ll 1+ \begin {cases} \log (B-A)/A &\text{si $\alpha=1$ ;} \\
            \frac{1}{1-\alpha} \Big( 1- \big(
\frac{B-A}A\big)^{\alpha-1}\Big) & \text{si $\alpha\neq 1$}
    \end{cases}  \\
&\ll \begin{cases} 1+\log (B/A) & \text{si $\alpha=1$ ;} \\
                1+ \big( \frac{B-A}A \big)^{ \alpha-1} &\text{sinon.} 
        \end{cases}
\end{align*} 
De plus, $\displaystyle \frac{1}{B-A}\leq \frac2B \leq \frac1A$,
si bien que
\begin{align*}
 I(A,B;\alpha,\beta) & \ll \frac{1}{(B-A)^{\alpha+\beta-1}}
     \times \left\{ \begin{matrix} 1+\log (B/A) \\
                            1+ ((B-A)/A)^{\alpha-1}\end{matrix}\right. \\
& \ll 
   \begin{cases} 
       (1+\log(B/A))/{A^{\alpha-1}B^\beta} &\text{si $\alpha=1$ ;} \\
        1/{A^{\alpha-1}B^\beta} & \text{sinon.}
    \end{cases}
\end{align*}
Le lemme est donc démontré.
\end {proof}

\begin{lemm}\label{lemm:-}%
Si $A,B\geq 1$, $\alpha\leq 1$, on a
\[ \int_0^{B-1} \frac{du}{(A+u)^\alpha (B-u)}
    \ll_\alpha \frac {1+\log A}{A^\alpha}. \]
\end{lemm}
\begin{proof}
On fait le changement de variables $A+u=(A+B)(1-t)$, soit
$B-u=(A+B)t$. Ainsi, l'intégrale vaut
\[ J(A,B;\alpha) =
\frac{1}{(A+B)^\alpha} \int_{1/(A+B)}^{B/(A+B)} \frac{du}{(1-u)^\alpha u} . 
\]
Si $A\geq B$, $u\leq B/(A+B)\leq 1/2$, donc $1-u\geq 1/2$ et
l'intégrale vérifie
\[ J(A,B;\alpha) \ll \frac{1}{(A+B)^\alpha} \int_{1/(A+B)}^{B/(A+B)}\frac{du}u
              = \frac{\log B}{(A+B)^\alpha}
        \ll \frac{1+\log A}{A^\alpha}. \]
Si $A\leq B$, on découpe l'intégrale de $1/(A+B)$ à $1/2$
et de $1/2$ à $B/(A+B)$. 
\begin{align*}
 \int_{1/(A+B)}^{1/2} \frac{du}{(1-u)^\alpha u}
& \leq \int_{1/(A+B)}^{1/2} \frac{du}u = \log\frac{A+B}2 \\
\int_{1/2}^{B/(A+B)} \frac{du}{(1-u)^\alpha u}
& \leq \begin{cases} \int_{1/2}^ 1(\cdots)  &\text{si $\alpha<1$ ;} \\
        \log \frac{A+B}{2A}\leq \log \frac{A+B}2 &\text{si $\alpha=1$}
        \end{cases}
\end{align*}
Finalement,
\[ J(A,B;\alpha)\ll \frac{1+\log(A+B)}{(A+B)^\alpha}\ll \frac{1+\log
A}{A^\alpha}, \]
ainsi qu'il fallait démontrer.
\end{proof}

\begin{proof}[Preuve du lemme~\ref{lemm:alpha}]
Si $a>0$, l'intégrale se majore par
\[ \int_0^\infty \frac{1}{(A+t)^\alpha}\frac{dt}{1+t}
\ll \frac{1+\log A}{A^\alpha} \]
d'après le lemme~\ref{lemm:+}.
Si $a<0$, on découpe l'intégrale de $0$ à $-a$ et de $-a$ à $+\infty$.
L'intégrale de $0$ à $-a$ vaut
\[
\int_0^{-a} \frac{1}{(A-t-a)^\alpha}\frac{dt}{1+t}
= \int_0^{-a}\frac{1}{(A+u)^\alpha}\frac{dy}{(1-a)-u}
\ll \frac{1+\log A}{A^\alpha}
\]
en vertu du lemme~\ref{lemm:-}, tandis que l'intégrale de $-a$ à $+\infty$
s'estime ainsi :
\[ \int_{-a}^\infty \frac{1}{(A+t+a))^\alpha}\frac{dt}{1+t}
=  \int_0^\infty \frac{1}{(A+u)^\alpha}\frac{du}{1-a+u}
\ll \frac{1+\log A}{A^\alpha}
\]
en appliquant de nouveau le lemme~\ref{lemm:+} et en distinguant suivant
que $A\leq 1-a$ ou $A\geq 1-a$.
\end{proof}

\egroup

\edef\guillemotleft{\guillemotleft\nobreak\,}
\edef\guillemotright{\unskip\nobreak\,\guillemotright}

\frenchspacing
\providecommand{\noopsort}[1]{}
\providecommand{\bysame}{\leavevmode ---\ }
\providecommand{\og}{``}
\providecommand{\fg}{''}
\providecommand{\smfandname}{et}
\providecommand{\smfedsname}{\'eds.}
\providecommand{\smfedname}{\'ed.}
\providecommand{\smfmastersthesisname}{M\'emoire}
\providecommand{\smfphdthesisname}{Th\`ese}

\end{document}